\documentclass[default,iicol]{sn-jnl}

\usepackage{graphicx}%
\usepackage{multirow}%
\usepackage{amsmath,amssymb,amsfonts}%
\usepackage{amsthm}%
\usepackage{mathrsfs}%
\usepackage[title]{appendix}%
\usepackage{xcolor}%
\usepackage{textcomp}%
\usepackage{manyfoot}%
\usepackage{booktabs}%
\usepackage{algorithm}%
\usepackage{algorithmicx}%
\usepackage{algpseudocode}%
\usepackage{listings}%
\usepackage{caption}
\usepackage{subcaption}
\usepackage{afterpage}
\usepackage{hyperref}
\usepackage{lscape}

\theoremstyle{thmstyleone}%
\theoremstyle{thmstyletwo}%

\theoremstyle{thmstylethree}%
\newcommand{\bm}{\boldsymbol}
\newcommand{\grad}{\bm{\nabla}}

\begin{document}

\title[Interpolation-based immersogeometric analysis methods for multi-material and multi-physics problems]{Interpolation-based immersogeometric analysis methods for multi-material and multi-physics problems}


\author[1]{\fnm{Jennifer E.} \sur{Fromm}}\email{jefromm@ucsd.edu}
\author[2]{\fnm{Nils} \sur{Wunsch}}\email{nils.wunsch@colorado.edu}
\author[2]{\fnm{Kurt} \sur{Maute}}\email{kurt.maute@colorado.edu}
\author[2]{\fnm{John A.} \sur{Evans}}\email{john.a.evans@colorado.edu}
\author*[3]{\fnm{Jiun-Shyan} \sur{Chen}}\email{js-chen@ucsd.edu}

\affil[1]{\orgdiv{Mechanical and Aerospace Engineering}, \orgname{University of California San Diego}, \orgaddress{\street{9500 Gilman Dr}, \city{La Jolla}, \postcode{92093}, \state{California}, \country{USA}}}

\affil[2]{\orgdiv{Aerospace Engineering}, \orgname{University of Colorado, Boulder}, \orgaddress{\street{3775 Discovery Drive}, \city{Boulder}, \postcode{80303}, \state{Colorado}, \country{USA}}}

\affil[3]{\orgdiv{Structural Engineering}, \orgname{University of California San Diego}, \orgaddress{\street{9500 Gilman Dr}, \city{La Jolla}, \postcode{92093}, \state{California}, \country{USA}}}

 
\abstract{Immersed boundary methods are high-order accurate computational tools used to model geometrically complex problems in computational mechanics.
While traditional finite element methods require the construction of high-quality boundary-fitted meshes, immersed boundary methods instead embed the computational domain in a structured background grid.
Interpolation-based immersed boundary methods augment existing finite element software to non-invasively implement immersed boundary capabilities through extraction.
Extraction interpolates the structured background basis as a linear combination of Lagrange polynomials defined on a foreground mesh, creating an interpolated basis that can be easily integrated by existing methods.
This work extends the interpolation-based immersed isogeometric method to multi-material and multi-physics problems.
Beginning from level-set descriptions of domain geometries, Heaviside enrichment is implemented to accommodate discontinuities in state variable fields across material interfaces.
Adaptive refinement with truncated hierarchically-refined B-splines is used to both improve interface geometry representations and to resolve large solution gradients near interfaces.
Multi-physics problems typically involve coupled fields where each field has unique discretization requirements.
This work presents a novel discretization method for coupled problems through the application of extraction, using a single foreground mesh for all fields.
Numerical examples illustrate optimal convergence rates for this method in both 2D and 3D, for partial differential equations representing heat conduction, linear elasticity, and a coupled thermo-mechanical problem. The utility of this method is demonstrated through image-based analysis of a composite sample, where in addition to circumventing typical meshing difficulties, this method reduces the required degrees of freedom when compared to classical boundary-fitted finite element methods. 
}






\maketitle

\section{Introduction}\label{sec:Intro}



Computational modeling has become an integral part of engineering of all types, and developments in manufacturing and design have increased the complexity of the models required. 
Multi-material problems are now ubiquitous, appearing in the design and analysis of composites \cite{rajak_fiber-reinforced_2019}, advanced additive manufacturing products \cite{nazir_multi-material_2023}, and multi-phase system analysis \cite{zhu_mixed_2021}. 
These new design spaces present challenges to modeling methods, primarily in the discretization of intricate geometries and domain interfaces, and the unique discontinuities in solutions that result from material interactions.
Finite element methods (FEMs) are among the most widely used computational tools for structural analysis. Classical FEM relies on sufficiently refined boundary-fitted meshes to discretize both the geometric domain of a material and its solution's function space. 
The accuracy of these methods is closely tied to mesh quality \cite{burkhart_finite_2013}, especially when high-order methods are employed \cite{engvall_mesh_2020}. Thus, considerable care must be taken in mesh construction before analysis can be performed \cite{knupp_algebraic_2001}. 
As the complexity of multi-material problems increases, generating these high-quality conforming meshes required by FEM becomes increasingly challenging, especially in three dimensions. 
Even for single material problems, mesh generation and refinement can consume up to 80\% of the design-through-analysis time of engineers \cite{boggs_dart_2005}. 
More complicated PDEs involving multiple state variables also present unique modeling challenges. In such multi-physics problems, the various fields may have distinct discretization needs and the fields must be coupled.
 
Immersed boundary methods circumvent conforming mesh generation by embedding the geometric problem domain into a background grid constructed on a geometrically simple domain. , Initially proposed to track fluid-structure interfaces in \cite{peskin_flow_1972}, similar classes of immersed methods were likewise developed by the solid mechanics community to accommodate discontinuities in solution fields without remeshing to create boundary-fitted meshes. The partition of unity method (PUM), introduced in \cite{babuska_partition_1997}, leverages the concept of enriching solution functions using \textit{a priori} knowledge of the location of discontinuities. This was combined with classical FEM in \cite{strouboulis_generalized_2000} and \cite{strouboulis_design_2000} to introduce the generalized finite element method (GFEM). A similar enriched method known as the "eXtended" finite element method (XFEM) \cite{moes_finite_1999, belytschko_elastic_1999, belytschko_arbitrary_2001} was also introduced to model crack propagation and other discontinuous problem. As opposed to the \textit{a priori} knowledge used to place enrichment functions used by GFEM, XFEM is characterized by adaptive enrichment schemes. Immersed boundary methods have also been extended to include high-order methods. The finite cell method, \cite{parvizian_finite_2007, schillinger_finite_2015} utilizes $p$ refinement of Lagrange polynomial basis functions to increase error convergence rates within an immersed framework. In the field of meshfree methods, the concepts of immersed or embedded methods have been used to model heterogenous materials \cite{schlinkman_quasi-conforming_2023}, and to enhance solution accuracy and stability near material interfaces in fluid-structure interaction problems \cite{huang_variational_2022}.

Isogeometric analysis (IGA) directly utilizes the geometric representation used in most computer-aided design (CAD) software in analysis \cite{hughes_isogeometric_2005} and was initially proposed to address the issue of generating high-quality boundary-fitted meshes. IGA has also been enhanced through the application of immersed boundary methods for the analysis of 'trimmed' CAD geometries \cite{burman_extension_2023}. The B-spline basis functions used in IGA offer additional advantages over classical FEM including improved geometric representation, higher levels of continuity, and improved per-degree-of-freedom accuracy compared to more common nodal finite element basis functions \cite{hughes_duality_2008, hughes_finite_2014}. These attributes make B-splines attractive choices for certain applications, including Kirchoff--Love shell analysis \cite{kiendl_isogeometric_2015}. The combinations of IGA and immersed boundary methods are often called immersogeometric methods \cite{kamensky_immersogeometric_2015}, and have been validated for use with hierarchically refined T-spline CAD models \cite{schillinger_isogeometric_2012}. The XFEM class of enriched methods was applied to IGA to create the eXtended isogeometric method (XIGA) in \cite{noel_xiga_2022}, exploiting level-set geometric descriptions and sophisticated integration algorithms to solve complicated PDEs involving multiple materials. XIGA was extended to use truncated hierarchically refined B-splines (THB-splines) for multi-physics problems in \cite{schmidt_extended_2023}. 

While offering elegant solutions to the problem of modeling complex geometries, existing immersed boundary software is currently limited to custom research codes. This is in part due to the complexity of generating custom quadrature rules on each cut background element \cite{divi_error-estimate-based_2020}. These custom quadrature rules require implementations that depart significantly from the highly optimized integration algorithms used in classical FEM codes. For example, an immersed boundary functionality called MultiMesh \cite{johansson_multimesh_2020} was developed for the popular open-source FEM software FEniCS \cite{alnaes_fenics_2015}, however it proved too difficult to maintain and was not ported to the more recent FEniCSx \cite{baratta_dolfinx_2023}. 

Classes of immersed boundary methods have been developed to reduce the difficulty of integrating over cut cells, including approximate domain methods and specifically the shifted boundary method \cite{main_shifted_2018, main_shifted_2018-1}. The shifted boundary method maps the boundaries of a computational domain to a mesh conforming surrogate domain and has been extended to high-order methods  in \cite{atallah_high-order_2022}.  Interpolation-based immersed boundary methods exist to address the same concerns. 

This work utilizes approximate extraction to retrofit classical FEM codes to perform immersed analysis.  
Extraction was introduced to represent B-spline basis functions with B\'ezier polynomials for implementation of IGA \cite{borden_isogeometric_2011}. 
B\'ezier extraction was generalized to Lagrange extraction in \cite{schillinger_lagrange_2016}. 
With Lagrange extraction, each spline in a B-spline function space is represented as a linear combination of Lagrange polynomials belonging to an interpolatory function basis $\{N_i\}$, which satisfies the Kronecker delta property $\delta_{ij} = N_i(\bm{x}_j)$, where $\bm{x}_j$ are the nodal coordinates. 
Because of this interpolatory property, Lagrange extraction-based methods are also called interpolation-based methods. 
IGA and other partition of unity methods are challenging to implement as shape functions are not the same from element to element \cite{moutsanidis_reduced_2021}> making interpolation-based methods attractive. 
Lagrange extraction has been utilized to implement IGA within existing finite element software in \cite{kamensky_tigar_2019, kamensky_open-source_2021}, which uses the software FEniCS, and in \cite{tirvaudey_non-invasive_2019}, which uses the software Code\_Aster. 

In \cite{fromm_interpolation-based_2023}, approximate extraction was applied to immersed boundary methods to address the same integration challenges presented by IGA. 
In interpolation-based immersed boundary methods, a domain is embedded in a structured background mesh, which is used to define a background basis. 
A boundary-fitted foreground mesh is quickly generated by decomposing the cut elements of the background. 
The foreground mesh, which does not need to adhere to the usual mesh quality metrics \cite{knupp_algebraic_2001}, is used to define a boundary-fitted Lagrange polynomial foreground basis. 
The background basis is interpolated with the foreground basis, and the interpolated basis can then be used to solve PDEs with Galerkin-type methods. 
As the interpolated basis is represented as a linear combination of Lagrange polynomials defined upon a boundary-fitted mesh, the integration can be performed with classical FE software, as illustrated with FEniCS in \cite{fromm_interpolation-based_2023}. 
This work also employed B-spline background bases, combining immersed and isogeometric methods for an immersogeometric method. 
The interpolation-based immersogeometric method differs from previous extraction-based IGA methods in its use of approximate instead of exact extraction. 
With approximate extraction, the interpolated background basis is not everywhere equivalent to the actual background basis, resulting in additional interpolation error. The additional interpolation error is demonstrated to be bounded by the optimal method error, thus the method still yields optimal convergence rates. 
Easing exact interpolation constraints to permit approximate extraction allows for easier and more efficient implementation.

This work expands interpolation-based immersogeometric methods to multi-material and multi-physics problems. In the previous work \cite{fromm_interpolation-based_2023}, a single uniform background basis was used for all fields of a single material problem. Here, a framework is presented to interpolate multiple enriched hierarchically refined background bases and extract these bases to a single foreground mesh. Hierarchical refinement allows for local refinement around material interfaces, which are described by level set functions. Heaviside enrichment applied at material interfaces allow the background bases to approximate discontinuous state variable fields. Within this new framework, separate background bases may be used to approximate different fields in multi-physics applications, which are interpolated using a single foreground basis, allowing for easy coupling. This framework is implemented in the next generation open-source software library FEniCSx \cite{baratta_dolfinx_2023}, with code available at \cite{fromm_jefrommexhume_dolfinx_2023}. 
Multi-materials heat conduction and linear elasticity are modeled to demonstrate the error convergence rates of the proposed workflow. A coupled thermo-mechanical problem illustrates the combined multi-material and multi-physics capabilities of the interpolation-based immersed boundary framework. 

The outline of this paper is as follows: Section \ref{sec:Bsplines} provides an overview of classic B-splines and the generation of the truncated hierarchically refined B-splines (THB-splines) used for local refinement. Section \ref{sec:MatInterfaces} describes this method's treatment of the geometric description of material interfaces and the Heaviside enrichment of solution spaces with discontinuities at material interfaces. Section \ref{sec:Interpolation} details the novel interpolation-based immersed boundary method's application to multi-material and multi-physics problems and its implementation workflow within existing FEM codes. Section \ref{sec:NumResults} will provide numerical results validating and expanding upon this method, and finally Section \ref{sec:Conclusions} will draw conclusions and indicate future work to be done with this method.  


\section{Hierarchical B-splines}\label{sec:Bsplines}
   Immersogeometric analysis combines aspects of two classes of methods, immersed methods and isogeometric methods. As in isogeometric methods, this work employs splines to represent state variable fields. The following is an overview of spline fundamentals, starting with a review of multivariate B-splines, proceeding with a discussion of hierarchical refinement, and concluding with the process for generating truncated hierarchically refined B-splines (THB-splines). 

\subsection{A review of multivariate B-spline function spaces}\label{subsec:Bsplines}

Recall first the construction of a 1D univariate B-spline basis with polynomial order $p$: $\left\{ B_{i,p}(\xi) \right\}_{i=1}^n$, where $n$ is the number of basis functions.  The domain is discretized with a knot vector $\Xi = \{ \xi_1, \xi_2,..., \xi_{n+p+1}\}$ such that $\left\{ \xi_i \right\}_{i=1}^{n+p+1} \subset \mathbb{R}$ and $\xi_1\leq \xi_2 \leq ...\leq \xi_{n+p+1}$. The functions are then constructed recursively from the piecewise constant basis function 
\begin{align}
    B_{i,0}(\xi) = \begin{cases}
        1, & \text{ if } \xi_i \leq \xi \leq \xi_{i + 1} \\
        0, & \text{ else}
    \end{cases} , 
\end{align}
using the Cox-de Boor recursion formula \cite{de_boor_calculating_1972}
\begin{align}
    \nonumber B_{i,p}(\xi) &= \dfrac{\xi - \xi_i}{\xi_{i + p} - \xi_i} B_{i,p-1}(\xi) \\
    & \quad + \dfrac{\xi_{i+p+1} - \xi}{\xi_{i+p+1} - \xi_{i+1}} B_{i+1,p-1}(\xi) .
\end{align}
If no interior knots are repeated, the basis is $C^{p-1}$-continuous at each knot in the interior domain, and $C^{\infty}$-continuous between the knots. The basis will form a partition of unity if the first and last  knots are repeated $p+1$ times. 
Higher dimension basis functions can be constructed by applying tensor-product operations to the univariate functions, such that 
\begin{align}
    B_{\bm{i},\bm{p}}(\bm{\xi}) = \prod_{m= 1}^{d_{p}} B_{i_m, p_m}^m (\xi^m),
\end{align}
where $d_p$ is the parametric space dimension, and there are $d_p$ knot vectors $\Xi^m = \{ \xi^m_1, \xi^m_2,..., \xi^m_{n_m+p_m+1}\}$, where $n_m$ is the number of basis functions and $p_m$ is the polynomial order of the $m^{th}$ parametric direction. Here $\bm{i} = \{i_1,..., i_{d_p}\}$ is a multi-index and $\bm{p} = \{p_1, ..., p_{d_p}\}$ is the vector of polynomial degrees. The B-spline basis of the set of these functions is denoted by $\mathcal{B}_{\bm{p}} :=\{B_{\bm{i},\bm{p}}\}$. 

\subsection{Hierarchically refined B-splines}\label{subsec:hb}

\begin{figure*}
     \centering
     \begin{subfigure}[b]{0.49\textwidth}
     \captionsetup{width=0.9\textwidth}
         \centering        \includegraphics[width=\textwidth]{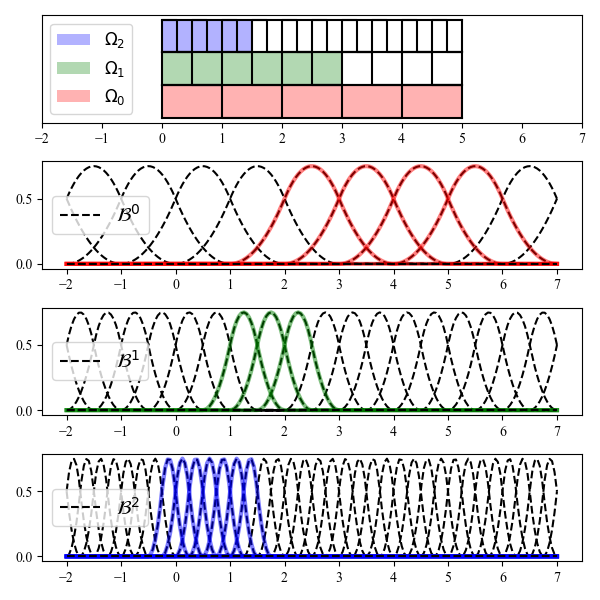}
         \caption{A 1D illustration of a sequence of subregions $\Omega_l$ and quadratic B-spline bases $\mathcal{B}^l$.}
         \label{fig:HB}
         \vspace{11pt}
     \end{subfigure}
     \begin{subfigure}[b]{0.49\textwidth}
     \captionsetup{width=0.9\textwidth}
         \centering
         \includegraphics[width=\textwidth]{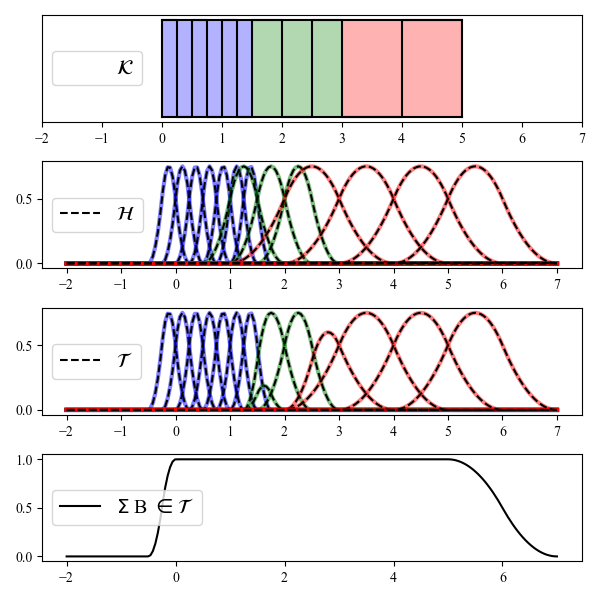}
         
         \caption{The hierarchically refined mesh $\mathcal{K}$, the HB-spline basis $\mathcal{H}$, the truncated basis $\mathcal{T}$, and the verification of the partition of unity.}
         \label{fig:THB}
     \end{subfigure}
     
        \caption{Local refinement is applied through truncated hierarchically-refined B-splines (THB). }
        \label{fig:HB_THB}
\end{figure*}
 Hierarchiacally refined B-splines (HB-splines) are constructed using nested sequences of spline spaces created by repeated knot insertion. Following the algorithms presented in \cite{garau_algorithms_2018}, an HB-spline basis begins with the construction of a sequence of $r$ tensor-product spline spaces $\mathcal{V}^l $,
\begin{align}
    \mathcal{V}^{0} \subset  \mathcal{V}^{1} \subset ... \subset \mathcal{V}^{r-1}, 
\end{align}
each of which has an accompanying B-spline basis $\mathcal{B}^l$, and tensor-product Cartesian mesh $\mathcal{K}^l$, where elements are denoted by $K$. A sequence of subdomains $\Omega^l$ are chosen, such that $\Omega^{l-1}$ is a subregion of $\Omega^l$, 
\begin{align}
    \Omega^{0} \supseteq \ \Omega^{1} \supseteq ... \supseteq \Omega^{r-1}, 
\end{align}
and each $\Omega^l$ can be discretized with mesh elements $K\in \mathcal{K}^l$. Here $r$ is the depth of refinement. 

The HB-spline basis $\mathcal{H} := \mathcal{H}^{r-1}$ is then constructed recursively by the algorithm:
\begin{align}
    \begin{cases}\label{eq:HAlg}
        \mathcal{H}^0 &:=\mathcal{B}^0,  \\
        \mathcal{H}^{l+1} &:= \{ B \in \mathcal{H}^l \ \big |  \ \text{supp}(B) \not\subset  \Omega_{l+1} \} \ \cup \\
        & \quad \ \{ B \in \mathcal{B}^{l+1} \ \big| \ \text{supp}(B)\subset  \Omega_{l+1} \}, \\
        & l \in \{ 0, ... , r-2\}.
    \end{cases}
\end{align}
In essence, each subsequent level's basis $\mathcal{H}^{l+1}$ is formed from the union of the set of basis functions from the previous level whose support is not in the new level's subdomain $\Omega^{l+1}$, and the set of functions from the new basis $\mathcal{B}^{l+1}$ whose support is within $\Omega^{l+1}$. This is illustrated with a 1D mesh in Figure \ref{fig:HB}.

The hierarchically refined basis $\mathcal{H}$ is associated with a hierarchically refined mesh $\mathcal{K}$, defined as
\begin{align}\label{eq:HRmesh}
    \mathcal{K} := \bigcup_{l=0}^{r-1} \{K \in \mathcal{K}^l \ \big |  \ K \in \Omega^l \text{ and } K \notin \Omega^{l+1} \}. 
\end{align}

\subsection{Enforcing the partition of unity property through truncation}\label{subsec:thb}

While a useful tool for applying adaptive refinement to IGA, HB-spline bases violate the partition of unity (PU) property. To regain the this property, the hierarchically refined bases are truncated as in \cite{garau_algorithms_2018} and \cite{giannelli_thb-splines_2012}. In addition to forming a partition of unity, truncation reduces the size of some basis functions' supports, thereby reducing the bandwidth of the resulting system of equations when compared to a non-truncated HB-spline basis.


A given multivariate basis function $B^l \in \mathcal{B}^l$ can be represented as a linear combination of the more refined functions of level $\mathcal{B}^{l+1}$:
\begin{equation}
    B^l = \sum_{B^{l+1} \in \mathcal{B}^{l+1}} c^{B^{l+1}}\big(B^l\big) \ B^{l+1}, 
\end{equation}
where $c^{B^{l+1}} \big(B^{l}\big)$ are coefficients relating the coarse basis function $B^l$ to the finer function $B^{l+1}$.

The truncation of $B^l$ removes the contributions from $B^{l+1} \in \mathcal{B}^{l+1}$ with support contained within $\Omega^{l+1}$, such that 
\begin{equation}
    \text{trun}^{l+1}(B^l) = B^l - \! \! \! \! \! \! \! \! \! \! \sum_{\substack{B^{l+1} \in \mathcal{B}^{l+1} , \\ \text{supp}(B^{l+1}) \subseteq \Omega^{l+1}}}  \! \! \! \! \! \! \! \! \! \! c^{B^{l+1}}\big(B^l\big) \ B^{l+1}.
\end{equation}

Using a similar algorithm to that given in Eq. \eqref{eq:HAlg}, the truncated basis $\mathcal{T} := \mathcal{T}^{r-1}$ can be constructed by 
\begin{align}
    \begin{cases}\label{eq:TAlg}
        \mathcal{T}^0 &:= \mathcal{B}^0 , \\
        \mathcal{T}^{l+1} &:= \{ \text{trun}^{l+1}(B)  \ \big |  \ B \in \mathcal{T}^l, \  \text{supp}(B) \not\subset  \Omega_{l+1} \} \\
        & \quad \cup \ \{ B \in \mathcal{B}^{l+1} \ \big |  \ \text{supp}(B)\subset  \Omega_{l+1} \}, \\
        & l \in \{ 0, ... , r-2\},
    \end{cases}
\end{align}
as illustrated in Figure \ref{fig:THB}.

\section{Immersed material interfaces} 
\label{sec:MatInterfaces}

Workflows to solve multi-material PDE require functionalities to both describe the geometry of material interfaces and to represent the associated discontinuities in the state variable fields. In this work, level set functions (LSFs) are utilized to implicitly describe the geometry of material interfaces, and a generalized Heaviside enrichment strategy in conjunction with a set of interface terms is employed to represent the required discontinuities at material interfaces.


\subsection{Representing interface geometry through level set functions} \label{subsec:LevelSet}

The level set method, developed in \cite{osher_fronts_1988}, has been used to describe interfaces in the extended finite element method (XFEM) \cite{moes_finite_1999, belytschko_elastic_1999}  and extended isogeometric analysis (XIGA) \cite{noel_xiga_2022, schmidt_extended_2023}. 
Following these works, the domain geometry is implicitly represented using LSFs $\phi_i(\bm{x})$. An iso-level $\phi_t$ of the LSF describes the interface $\Gamma_\pm$ between two subdomains $\Omega_+$ and $\Omega_-$ such that
\begin{align}
        \nonumber \phi(\bm{x}) < \phi_t, \text{  }& \bm{x} \in \Omega_{+}, \\
        \phi(\bm{x}) >\phi_t, \text{  }& \bm{x} \in \Omega_{-} ,\\
        \nonumber \phi(\bm{x}) = \phi_t, \text{  }& \bm{x} \in \Gamma_{\pm}.
\end{align}
With $n$ LSFs, this method can represent up to $2^n$ subdomains. Materials are then associated with these subdomains using  a multi-phase level set model as in \cite{vese_multiphase_2002}, where phases are identified by phase indices $\mathcal{P}$. Phase indices $\mathcal{P}$ are assigned with characteristic functions $f_i$, 
\begin{equation}
    f_i(\bm{x}) = \begin{cases}
        0,  \text{  }& \phi_i(\bm{x}) < \phi_t, \\
        1,\text{  }& \phi_i(\bm{x}) \geq \phi_t,\\
    \end{cases}
\end{equation}
such that 
\begin{equation}
\label{eq:phase_IDs}
    \mathcal{P}(\bm{x}) = \sum_{j=1}^n 2^{j-1} f_j(\bm{x}). 
\end{equation}
Phases are then mapped onto material subregions. 

In this work LSFs are discretized using linear basis functions from a THB-spline basis $B_k \in \mathcal{T}$,
\begin{equation}
    \phi^h_i (\bm{x}) = \sum_{k} B_k(\bm{x}) \phi_i^k ,
\end{equation}
where $\phi_i^k$ are the coefficients associated with LSF $\phi_i$. The LSF are linearly interpolated such that the coefficients are the nodal values $\phi^i_k= \phi^i(\bm{x}_k)$. 
This discretization is used to construct material indicator functions and to enrich background basis functions. 

\subsection{Heaviside enrichment of basis functions at material interfaces} 
\label{subsec:Heaviside}
\begin{figure}[b]
\centering
    \includegraphics[width=0.9\linewidth]{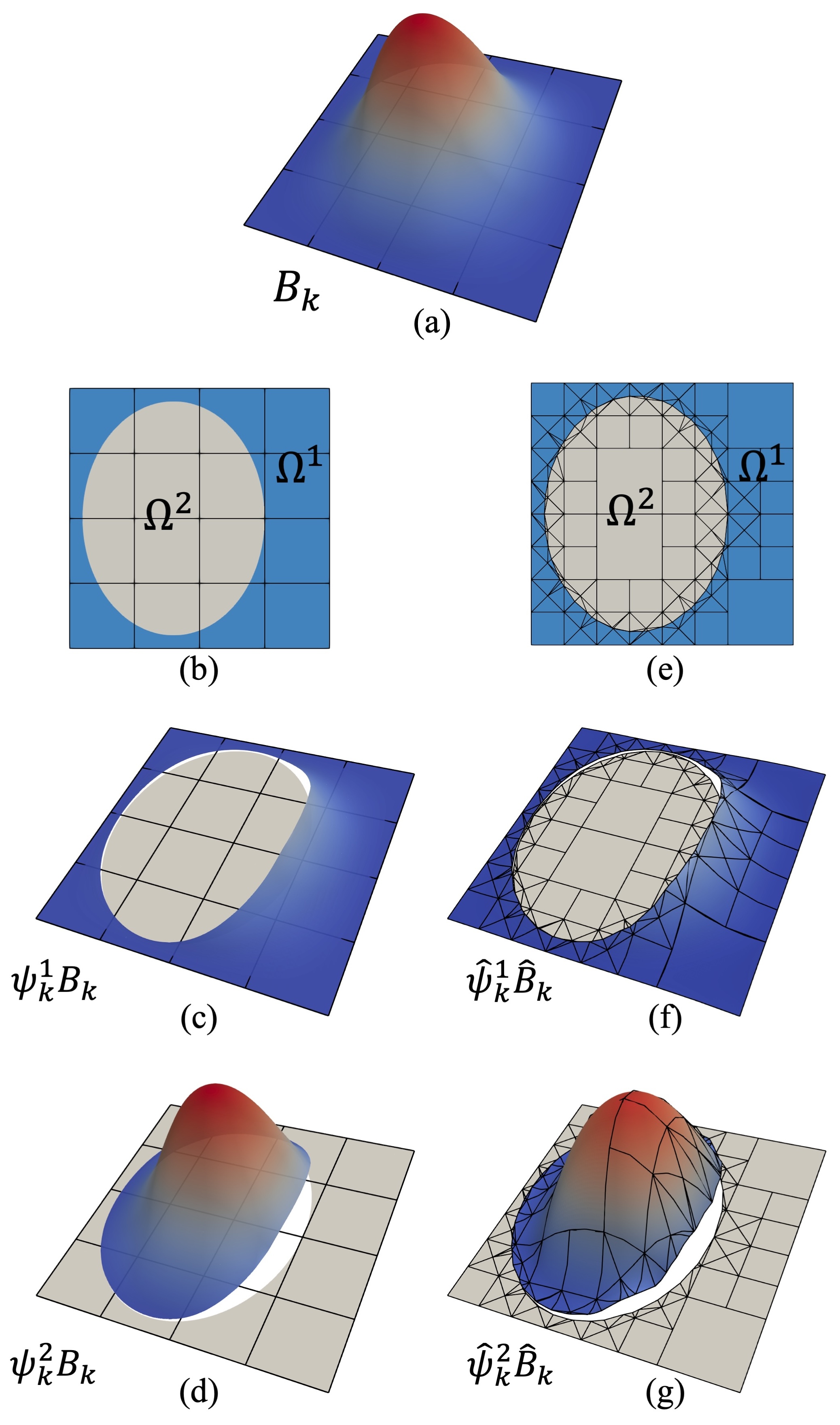}
    \caption{A B-spline basis is defined on a structured background mesh and an example function $B_k$ is depicted in (a). Using the geometry description of the material subdomains in (b), Heaviside enrichment is applied to form the discontinuous functions $\psi^1_k B_k$ and $\psi^2_kB_k$, depicted in (c) and (d) respectively. A Lagrange foreground function space is defined on the boundary-fitted mesh in (e). The function space is used to interpolate the enriched background functions $\Hat{\psi}^1_k \Hat{B_k}$ and  $\Hat{\psi}^2_k \Hat{B_k}$, depicted in (f) and (g), respectively. }
        \label{fig:enrich_int_ref}
\end{figure}
Heaviside enrichment have been widely used in PUFEM \cite{babuska_partition_1997}, GFEM \cite{strouboulis_generalized_2000}, and XFEM \cite{belytschko_arbitrary_2001} as a means to represent strong discontinuities within elements. The enrichment strategies presented in most existing literature, such as in \cite{terada_finite_2003, hansbo_finite_2004}, 
adds enriched basis functions for each material.
This work adopts the enrichment strategy presented by \cite{noel_xiga_2022} which instead considers the material connectivity in the individual basis functions' supports to alleviate artificial numerical stiffening arising around small geometric features. 
This stiffening is caused by interpolating the state variable field in locally disconnected domains of the same material by the same basis function. 
The high-order, higher-continuity B-spline basis functions with large supports employed in this work, alongside the complex material layouts presented in Subsection \ref{subsec:composite_image}, would exacerbate local stiffening effects and lead to an increased solution error if a global enrichment strategy were instead chosen. Additionally, function wise enrichment results in fewer added degrees of freedom than global enrichment, which reduces overall system size. 


For a given function $B_k$ with support $\textup{supp}(B_k)$, the phase IDs, defined in Eq. \eqref{eq:phase_IDs}, are used to identify the $L_k$ distinct but connected material subregions $\Omega_k^m$, such that supp$(B_k) = \cup_{m=1}^{L_k} \Omega_k^m$. For $L_k$ distinct subregions, the basis function $B_k$ requires $L_k$ enrichment levels. 

This enrichment is then achieved through indicator functions $\psi_k^m$, 
\begin{equation}
    \psi_k^m(\bm{x}) = \begin{cases}
        1, \text{ if } \bm{x} \in \Omega_k^m\\
        0, \text{ else,}
    \end{cases}
\end{equation}
such that the enriched basis functions can be expressed as 
\begin{align}\label{eq:enrichB}
    B^m_k(\bm{x}) = \psi_k^m (\bm{x}) B_k(\bm{x}), \ \ \forall \ m \in \{1, ..., L_k\}. 
\end{align}

The enriched basis functions constructed from a single non-enriched basis function, are shown in Figure \ref{fig:enrich_int_ref} for a two-material configuration. The bi-quadratic B-spline $B_k$ depicted in Figure \ref{fig:enrich_int_ref}(a) is enriched assigning one material to the inside and one material to the outside of the ellipse shown in
Figure \ref{fig:enrich_int_ref}(b). The basis function $B_k$ is split into two enriched functions, $B^1_k(\bm{x}) = \psi_1^m (\bm{x}) B_k(\bm{x})$  in Figure \ref{fig:enrich_int_ref}(c) and $B^2_k(\bm{x}) = \psi_k^2 (\bm{x}) B_k(\bm{x})$ in Figure \ref{fig:enrich_int_ref}(d). allowing for the representation of strong discontinuities at the material interface. Interface conditions are enforced weakly; for example $C^0$ continuity can be enforced at the interface using Nitsche's method \cite{annavarapu_robust_2012}. 

\section{The interpolation-based immersed boundary method}
\label{sec:Interpolation}
One of the core challenges associated with classic immersed methods is the construction of custom quadrature rules for the various material regions within each intersected background element.
Numerous solutions exist and have been used in custom research codes, such as octree refinement \cite{duster_finite_2008, schillinger_finite_2015, kamensky_immersogeometric_2015}, interface reconstruction and tessellation \cite{min_geometric_2007,cheng_higher_2010,fries_higher_2016}, moment-fitting \cite{muller_highly_2013, sudhakar_quadrature_2013}, and quadrature schemes using generalized Stokes theorem \cite{sommariva_gauss_2009,gunderman_spectral_2021,gunderman_high_2021}. However, the generation of such quadrature rules can generally not be implemented within existing finite element software without major changes to the software itself. 

The interpolation-based immersed approach is instead designed to utilize the integration subroutines of existing FEM codes. 
To this end, a boundary-fitted foreground mesh is constructed with only minimal requirements on mesh quality. Element formation is then performed on the poor quality boundary-fitted mesh using existing standard finite element routines. Using Lagrange extraction operators \cite{schillinger_lagrange_2016} the resulting tangent matrix and force vector are projected into the enriched THB-spline space. This can be done either either on an elemental level during assembly, or globally afterwards. 
The resulting final problem uses an approximation of the enriched function space of the background mesh which is interpolated by the basis functions of the foreground mesh. 

The following Subsection \ref{subsec:ModelProblem} will first introduce the thermo-elastic model problem. Subsection \ref{subsec:InterpolatedBF} will provide an overview of the interpolation-based approach specific to the multi-material, multi-physics problems presented in this paper. For a more general and comprehensive introduction to the approach, we refer the reader to the authors' previous work on the topic \cite{fromm_interpolation-based_2023}. The generation of the boundary-fitted foreground mesh is discussed in Subsection \ref{subsec:FGMesh}.



\subsection{Multi-material and multi-physics model problem} 
\label{subsec:ModelProblem}

To illustrate the application of interpolation-based immersed boundary methods to multi-material and multi-physics problems, a thermo-elastic problem is introduced. This problem can be broken into a thermal subproblem and a structural subproblem. 

Let a domain of interest $\Omega$ with closure denoted $\overline{\Omega}$ be composed of $n$ material subdomains $\mathcal{M} = \{1, \cdots, n \}$
\begin{align}
    \overline{\Omega} = \bigcup_{m \in \mathcal{M} } \overline{\Omega}^{m} \quad \subset \mathbb{R}^d.
\end{align} 
A different thermal conductivity $\kappa^m$ may be associated with each material $m$
\begin{align}
\kappa(\bm{x}) = \kappa^m \text{, } ~ & \bm{x} \in \Omega^m,
\end{align}
for $m \in \mathcal{M}$. A source term $f:\Omega\rightarrow\mathbb{R}$, a boundary heat flux term $\overline{q}:\partial \overline{\Omega}\rightarrow\mathbb{R}$ on $\Gamma_{\Bar{q}} \!\subset \! \partial \overline{\Omega}$, and Dirichlet boundary data $\overline{T}:\partial \overline{\Omega}\rightarrow\mathbb{R}$ on $\Gamma_{\overline{T}} \! \subset \!\partial \overline{\Omega}$ are ascribed. The strong form for the thermal problem then reads as:  Find $T: \Omega  \rightarrow \mathbb{R}$ such that $\forall$ $m\in \mathcal{M}$
\begin{equation}
\begin{split}
\label{eq:poisson-strong}
    -\grad \bm{\cdot} ( \kappa(\bm{x})\grad T ) = f   ~~~ & \text{ in } \Omega^m \text{,} \\
    [\![ T]\!] = 0   ~~~ &  \text{ on all } \Gamma_{km}, \\
    [\![ \bm{q}]\!] = 0   ~~~ &  \text{ on all } \Gamma_{km}, \\
     - \kappa(\bm{x}) \grad T \cdot \bm{n}  =\overline{q} ~~~ & \text{ on } \Gamma_{\Bar{q}}^m \text{,} \\        
     T  = \overline{T}  ~~~ & \text{ on } \Gamma_{\overline{T}}^m \text{,}\\
\end{split}
\end{equation}
where $\Gamma_{km }= \overline{\Omega}^k \cap \overline{\Omega}^m \neq \emptyset$, with $k\in\mathcal{M}$ and $ k \neq m$, are the material interfaces, and $[\![ \cdot]\!] = (\cdot)^{k} - (\cdot)^{m}$ is the jump of a given quantity over an interface $\Gamma_{km}$.  The material fields are defined $T^m =T(\bm{x})$, $\bm{x} \in \Omega^m$, and  $ \bm{q}^m =  - \kappa^m \grad T^m$. 
The domains $\Gamma_{\Bar{q}}^m = \Gamma_{\Bar{q}} \cup \partial \overline{\Omega}^m$, and $\Gamma_{\overline{T}}^m = \Gamma_{\overline{T}} \cup \partial \overline{\Omega}^m$ are the intersections of the domain boundaries with the material subdomain boundaries. 
$\bm{n}$ denotes the surface normal. 

The domain of interest $\overline{\Omega}$ is embedded into a hierarchically refined background mesh $\mathcal{K}_T$, generated using the sequence of refined meshes $\mathcal{K}^l$ and subdomains $\Omega_{T}^l$.\footnote{Here the subscript $(\cdot)_{T}$ refers to entities associated with the thermal subproblem.  The subscript $(\cdot)_{u}$ will refer to entities associated with the structural subproblem. $T$ and $u$ are used as superscripts when referring to entities that require subscripts for indices, such as  B-spline basis functions ($B^T_i$ and $B^u_i$).} Note that the largest subdomain $\Omega_{T}^0$ must be chosen such that the closure of the domain of interest is a subset, $\overline{\Omega}\subset \Omega_T^0$. The mesh $\mathcal{K}_T$ is associated with the enriched THB-spline basis $\mathcal{T}_T = \{B^T_{i}\}$. The temperature field is discretized using the function space 
\begin{align}
    \mathcal{V}^h_T = \text{span}\{B^T_i \ \big| \ \text{supp}(B^T_i) \cap \overline{\Omega}\neq \emptyset\}.
\end{align}

The discrete form can then be defined as: Find $T^h \in \mathcal{V}^h_T$ such that $\forall$ $\theta^h\in\mathcal{V}^h_T$,
\begin{align}
& \sum_{m=1}^n \left[\int_{\Omega^m} \! \! \kappa \grad T^h \bm{\cdot} \grad \theta^h d\Omega\right]   - \int_{\Omega} f \theta^h d\Omega  -  \int_{\Gamma_{\overline{q}}} \overline{q} \theta^h d\Gamma \nonumber \\
& \quad \quad = \mathcal{R}^{D}_T + \mathcal{R}^I_T \text{  ,}\label{eq:heat-disc-generic}
\end{align}
where $\mathcal{R}^{D}_T$ and $ \mathcal{R}^{I}_T$ are Dirichlet and interface residual terms. The temperature Dirichlet residual is the result of Nitsche's method \cite{nitsche_uber_1971} enforcement of the Dirichlet boundary condition, 
\begin{align} \label{eq:temp-nitsches}
    \nonumber \mathcal{R}^{D}_T &=   \sum_{m=1}^n \Bigg [\mp \int_{\Gamma_{\overline{T}}^m} \kappa (T^h - \overline{T}) (\grad \theta^h \cdot\bm{n}) \,d\Gamma  \\
    \nonumber &\quad -\int_{\Gamma_{\overline{T}}^m} \kappa \theta^h(\grad T ^h\cdot\bm{n})\,d\Gamma \\
    &\quad + \int_{\Gamma_{\overline{T}}^m}\frac{\beta^{D}_T \kappa }{h}(T^h - \overline{T} ) \theta^h\,d\Gamma \Bigg], 
\end{align}
where $\beta^{D}_T \geq 0$ is a user defined constant. The first integral of Eq. \eqref{eq:temp-nitsches} will be negative for the symmetric version of Nitsche's method (which is employed in this work) or positive for the non-symmetric version.  $h$ is taken as the characteristic element size on the foreground mesh, differing from the usual implementation where $h$ is the element size on the background mesh. 

The temperature interface conditions lines 2 and 3 of Eq. \eqref{eq:poisson-strong},  are also enforced through a Nitsche-like method, resulting in the temperature interface residual
\begin{align}
    \nonumber \mathcal{R}^{I}_T &=  \sum\limits_{i=1}^n \sum\limits_{j=i+1}^n \Bigg[- \int_{\Gamma_{ij}} [\![ T^h ]\!] \{ \kappa \grad \theta^h \} \cdot\bm{n}) \,d\Gamma \\
    \nonumber &\quad - \int_{\Gamma_{ij}} [\![ \theta^h ]\!] \{\kappa \grad T ^h\} \cdot\bm{n})\,d\Gamma \\
    &\quad + \int_{\Gamma_{ij}} \gamma^{ij}_T  [\![ T^h ]\!] [\![ \theta^h]\!] \,d\Gamma\Bigg],
\end{align}
where $\{\cdot\} = w^{i}(\cdot)^{i} - w^{j}(\cdot)^{j}$ is the weighted average of a given quantity. Motivated by the formulation in \cite{annavarapu_robust_2012}, these weights are defined as 
\begin{align} \label{eq:weights}
   \nonumber  w^{i} &= \dfrac{ (h^{i})^{d_p} /  \omega^{i} }{(h^{i})^{d_p}  / \omega^{i} + (h^{j})^{d_p} / \omega^{j} } \ \text{ and } \\
    w^{j} &= \dfrac{ (h^{j})^{d_p} /  \omega^{j} }{(h^{i})^{d_p}  / \omega^{i} + (h^{j})^{d_p} / \omega^{j} } , 
\end{align}
where $h^{m}$ is the characteristic size of the foreground element in domain $\Omega^{m}$ bordering the interface facet, $\omega^m$ is the characteristic material parameter, which for the thermal subproblem is $\kappa^m$, and $d_p$ is the domain dimension. The penalty parameter $\gamma^{ij}_T  $ is defined as 
\begin{equation}
    \gamma^{ij}_T  = 2 \beta^{I}_T  \dfrac{(h^{i})^{d_p-1} + (h^{j})^{d_p-1}}{(h^{i})^{d_p}  / \omega^{i} + (h^{j})^{d_p} / \omega^{j} }, 
\end{equation}
 where $\beta^{I}_T \geq 0$ is a user specified constant.

 The thermal subproblem can be stated compactly as the variational problem: Find $T^h \in \mathcal{V}^h_T$ such that $\forall \ \theta^h\in\mathcal{V}^h_T$
 \begin{align}\label{eq:varT}
     a_T(T^h, \theta^h) = L_T(\theta^h),
 \end{align}
where $a_T(T,\theta)$ and $L_T(\theta)$ can be computed from Eq. \eqref{eq:heat-disc-generic}.

The structural subproblem may utilize a differently refined background discretization.  In this case a separate sequence of refined domains $\Omega_{u}^l$ may be selected. The structural subproblem's background mesh $\mathcal{K}_u$ is constructed withthis sequence of subdomains and the same sequence of refined tensor-product Cartesian grids $\mathcal{K}^l$. The basis $\mathcal{T}_u = \{B^u_{i}\}$ associated with the mesh is used for the components of the displacement field $\bm{u}$. Each displacement component is discretized using the function space 
\begin{align}
    \mathcal{V}^h_u = \text{span}\{B^u_i\ \big| \ \text{supp}(B^u_i) \cap \Omega \neq \emptyset\}. 
\end{align}

Using a similar derivation (given in full in Appendix \ref{app:MMlinearElasticity}) as applied to the temperature subproblem the mechanical variation problem is compactly written as: Find $\bm{u}^h \in \bm{\mathcal{V}}^h_u = [\mathcal{V}^h_u, \mathcal{V}^h_u]$ such that, $\forall  \bm{v}^h\in \bm{\mathcal{V}}^h_u$
 \begin{align}\label{varU}
     a_u(\bm{u}^h,\bm{v}^h) = L_u(\bm{v}^h).
 \end{align}

The two subproblems are coupled through a constitutive model accounting for the thermal expansion by computing the mechanical strain $\bm{\varepsilon}_{m}$ as 
\begin{align}\label{eq:TEstrain}
    \nonumber \bm{\varepsilon}_m (\bm{u},T) &= \bm{\varepsilon}_{u}(\bm{u}) - \bm{\varepsilon}_T(T) \\
    &= \dfrac{1}{2}\left( \grad \bm{u}+  (\grad\bm{u})^{\text{T}}  \right) - \alpha (T - T_0)\bm{I}, 
\end{align}
where $\alpha$ is the thermal expansion coefficient and $T_0$ is the temperature in the reference configuration. $\bm{I}$ is the identity matrix.

The fully coupled system is then given by the variational problem: Find $(\bm{u}^h,T^h)  \in [\bm{\mathcal{V}}^h_u, \mathcal{V}^h_T]$ such that, $\forall$ $(\bm{v}^h,\theta^h)\in [\bm{\mathcal{V}}^h_u, \mathcal{V}^h_T]$
 \begin{align} \label{eq:varUT}
    \nonumber a_T(T^h, \theta^h) &= L_T(\theta^h) \text{ and } \\
     a_u(\bm{u}^h, \bm{v}^h) - b (T^h,\bm{v}^h) &= L_u(\bm{v}^h),
 \end{align}
where the form $b (T,\bm{v})$ is a result of the coupling conditions and is given by 
\begin{align} \label{eq:bCoupled}
     \nonumber b(T,\bm{v}) &= \sum_{m=1}^n \Bigg [ \int_{\Omega^m}\bm{\varepsilon}^u(\bm{v}) : \bm{\varepsilon}^T( T) \text{d}\Omega \Bigg] \\
     \nonumber &+ \int_{\Gamma_{\overline{u}}} \bm{C}:\bm{\varepsilon}^T(T) \cdot \bm{n} \cdot \bm{v}^h \text{d}\Gamma \\
     &-\sum\limits_{i=1}^n \sum\limits_{j=i+1}^n \Bigg[ \int_{\Gamma_{ij}} [\![ \bm{v} ]\!] \cdot \{\bm{C}:\bm{\varepsilon}^T(T)\} \cdot\bm{n}\,d\Gamma \Bigg].
\end{align}

\subsection{Interpolated basis functions} 
\label{subsec:InterpolatedBF}



In traditional immersed boundary methods custom quadrature rules would be used to evaluate the integral in the weak form of the coupled problem, given in Eq. \eqref{eq:varUT}. The main idea of the interpolation-based immersed paradigm is to interpolate the background basis functions using a space of Lagrange functions defined on a foreground mesh, which can be integrated with classical quadrature methods. This workflow thus introduces an interpolated background function space for the thermal subproblem 
\begin{align}
    \label{eq:interpolated_temp_space}
    \widehat{\mathcal{V}}^h_T &= \text{span}\{\widehat{B}^T_i \ \big| \ \text{supp}(\widehat{B}^T_i) \cap \Omega \neq \emptyset\}, 
\end{align}
where the interpolated basis functions are defined as 
\begin{align} 
    \widehat{B}^T_i &:=\sum_{j=1}^\nu M_{ij}^T N_j
\end{align}
where 
\begin{align}
    M^T_{ij} := B^T_i(\bm{x}_j)
\end{align}
is the Lagrange extraction operator. $\{N_j\}_{j=1}^\nu$ is the basis of a Lagrange FE space with nodal points $\bm{x}_j$ such that $N_i(\bm{x}_j) = \delta_{ij}$. Here $\nu$ is the number of foreground basis functions. The same foreground space is used to interpolate the background bases for both the temperature and the displacement state variables. 

The approximations of the temperature field becomes 
\begin{align} 
    \label{eq:interpolationT}
    T^h &= \sum_{i=1}^{n_T} \widehat{B}^T_i  d^T_i= \sum_{i=1}^{n_T} \sum_{j=1}^{\nu} M^T_{ij} N_j d^T_i ,
\end{align}
where $\{d^{T}_i\}_{i=1}^{n_{T}}$ are the unknown coefficients and $n_T$ is the number of basis functions in the temperature field's interpolated background B-spline basis. The displacement is similarly discretized with vector valued function spaces as 
\begin{align}\label{eq:interpolationu_main}
    \bm{u}^h  &= \sum_{I=1}^{(d_p\cdot n_u)} \widehat{\bm{B}}^u_I  d^{\bm{u}}_I= \sum_{I=1}^{(d_p\cdot n_u)} \sum_{J=1}^{(d_p\cdot \nu)} M^{\bm{u}}_{IJ} \bm{N}_J d^{\bm{u}}_I.
\end{align}
Details regarding this discretization are given in Appendix \ref{app:MMlinearElasticity}.

The variational problem in Eq. \eqref{eq:varUT} can be assembled using the interpolated bases in Eqs. \eqref{eq:interpolationT} and \eqref{eq:interpolationu_main} to form the linear system 
\begin{align}
    \begin{bmatrix}
        \bm{K}^{\theta \theta} & \bm{0} \\
        \bm{K}^{\theta \bm{v}} &  \bm{K}^{\bm{v} \bm{v}}
    \end{bmatrix} \begin{bmatrix}
        \bm{d}^T \\
        \bm{d}^{\bm{u}}
    \end{bmatrix} = 
    \begin{bmatrix}
        \bm{f}^{\theta} \\
        \bm{f}^{\bm{v}}
    \end{bmatrix}, 
\end{align}
where 
\begin{align}
    K^{\theta \theta}_{ij}  &= a_{T}(\widehat{B}^T_i, \widehat{B}^T_j), \\
    K^{\theta \bm{v}}_{iJ}  &= b(\widehat{B}^T_i, \widehat{\bm{B}}^{\bm{u}}_J), \\
    K^{\bm{v} \bm{v}}_{IJ}  &= a_{\bm{u}}(\widehat{\bm{B}}^{\bm{u}}_I, \widehat{\bm{B}}^{\bm{u}}_J), \\
    f_i^{\theta} &= L_T( \widehat{B}^T_i), \quad \text{ and }\\
    f_i^{\bm{v}} &= L_u(\widehat{\bm{B}}^{\bm{u}}_I)
\end{align}
are the matrix entries. 

By consolidating the extraction operators for all components in a single matrix
\begin{align}
    \bm{M} = \begin{bmatrix}
        \bm{M}^{T} & \bm{0} \\
        \bm{0} & \bm{M}^{\bm{u}}
    \end{bmatrix}
\end{align}
the linear system can be rewritten as 
\begin{align}
    \bm{M}^{\text{T}}
    \begin{bmatrix}\label{eq:LinSystem}
        \bm{A}^{\theta \theta}& \bm{0} \\
        \bm{A}^{\theta \bm{v}} &  \bm{A}^{\bm{v} \bm{v}}
    \end{bmatrix} \bm{M}
    \begin{bmatrix}
        \bm{d}^T \\
        \bm{d}^{\bm{u}}
    \end{bmatrix} = 
    \bm{M}^{\text{T}}
    \begin{bmatrix}
        \bm{b}^{\theta} \\
        \bm{b}^{\bm{v}}
    \end{bmatrix}, 
\end{align}
where 
\begin{align}
\label{eq:global_system_with_extraction_operators}
    A^{\theta \theta}_{ij}  &= a_{T}(N_i, N_j), \\
    A^{\theta \bm{v}}_{iJ}  &= b(N_i, \bm{N}_J), \\
    A^{\bm{v} \bm{v}}_{IJ}  &= a_{\bm{u}}(\bm{N}_I, \bm{N}_J), \\
    b_i^{\theta} &= L_T(N_i), \quad \text{ and }\\
    b_I^{\bm{v}} &= L_u(\bm{N}_I).
\end{align} 

The quantities $\bm{A}^{\theta \theta}$, $\bm{A}^{\theta \bm{v}}$, $\bm{A}^{\bm{v} \bm{v}}$, $\bm{b}^{\theta}$, and $\bm{b}^{\bm{v}}$ are evaluated and assembled on the boundary-fitted foreground mesh. As the foreground basis is thge typical conforming Lagrange polynomial basis, existing commercial or open-source FE software applying standard quadrature rules may be utilized.

The extraction operators can either be computed globally, as suggested in Eq. \eqref{eq:global_system_with_extraction_operators}, or on an element level for more efficient implementation. However, the latter option requires a modification of the assembly maps in existing FE software. After solving the linear system in Eq. \eqref{eq:LinSystem} with the interpolated basis, the solution is post-processed by projecting the solution for each field onto the foreground basis. 

A visual example of the interpolation of an enriched B-spline basis function is given in Figure \ref{fig:enrich_int_ref}. The boundary-fitted foreground mesh in Figure \ref{fig:enrich_int_ref}(e) is used to define a discontinuous Lagrange polynomial foreground function space. The enriched background functions $B_k^1$ and $B_k^2$ are interpolated with this function space, as shown in Figure \ref{fig:enrich_int_ref}(f) and (g). 

\subsection{Foreground mesh generation}
\label{subsec:FGMesh}

\begin{figure*}
    \centering
    \includegraphics[width=0.9\linewidth]{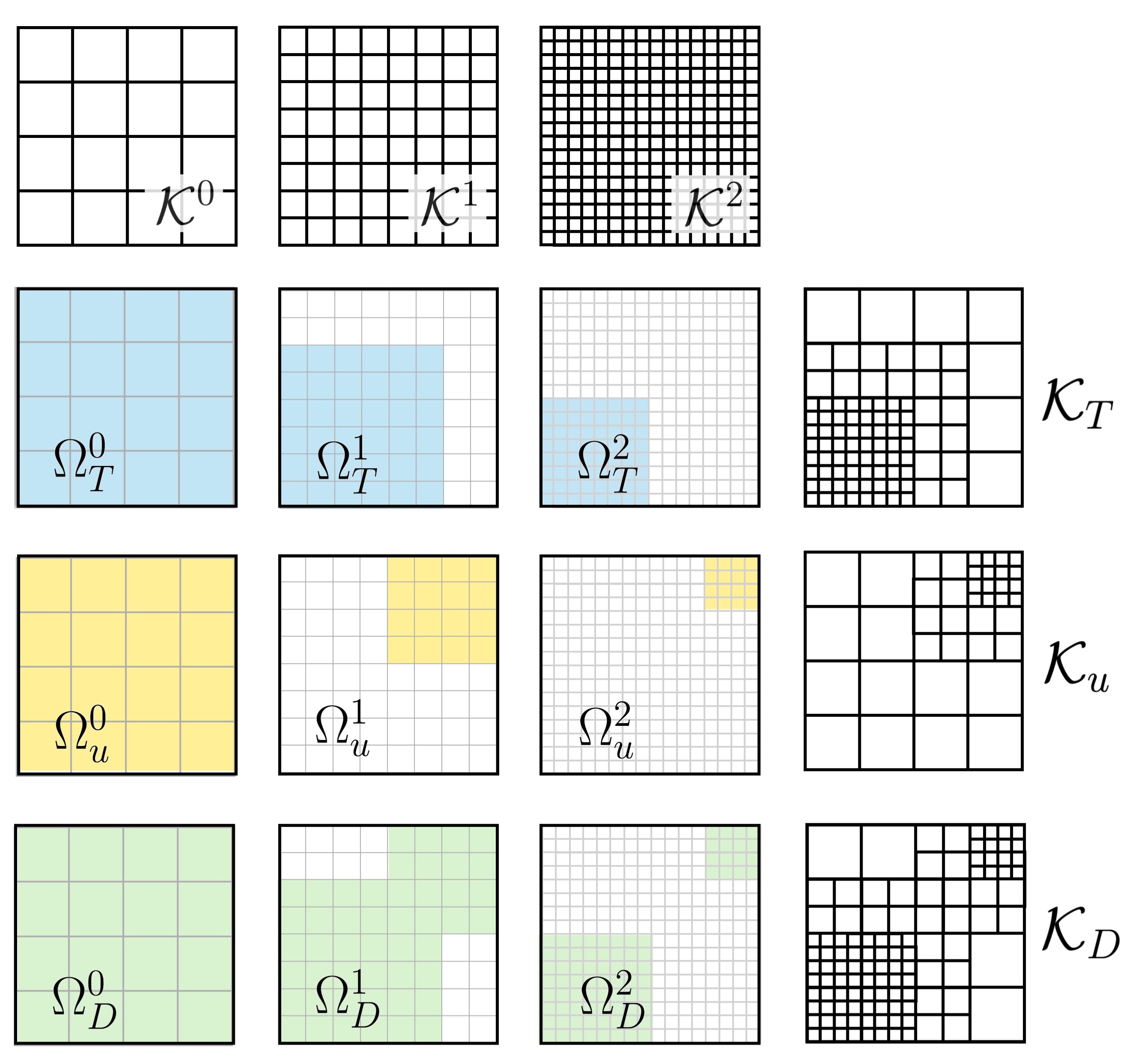}
    \caption{A series of uniformly refined meshes $\mathcal{K}^l$ are shown in the top row. A series of nested subdomains, $\Omega^{l+1}_T \subseteq \Omega^l_T$, shown in the second row, and $\Omega^{l+1}_u \subseteq \Omega^l_u$, shown in the third row, are defined for each state variable. For the decomposition mesh the series of subdomains $\Omega^l_{D}$, shown in the last row, is defined such that the domain on each level $l$ contains the union of the $l^{\text{th}}$  level domains for both each variable field. The hierarchically refined meshes $\mathcal{K}_T$, $\mathcal{K}_u$, and $\mathcal{K}_D$, shown in the rightmost column, are constructed with the mesh series $\mathcal{K}^l$ in the top row and their respective subdomain series.  }
    \label{fig:HRmesh}
\end{figure*}
Both state variable fields are interpolated using the same Lagrange FE space $\{N_j\}_{j=1}^\nu$ which is constructed on a boundary-fitted foreground mesh. To generate this foreground mesh, the elements of a background mesh intersected by interfaces are decomposed into triangles and tetrahedrons whose facets approximately reconstruct the interfaces. 

To ensure the foreground mesh is sufficiently refined, a background mesh for decomposition $\mathcal{K}_{D}$ is constructed from the series of hierarchically refined meshes $\mathcal{K}^l$,  
\begin{align}
    \label{eq:HRmesh_union}
    \mathcal{K}_{D} := \bigcup_{l=0}^{r-1} \{K \in \mathcal{K}^l| K \in \Omega_{D}^l \text{ and } K \notin \Omega_{D}^{l+1} \}, 
\end{align}
where $\Omega_{D}^l \supseteq \Omega_{T}^l \cup \Omega_{u}^l $ is a series of subdomains containing the union of the subdomains used for each subproblem's discretization. This ensures that $\mathcal{K}_{D}$ is at least as refined as the meshes used for each subproblem's discretization and allows for additional refinement of the foreground to improve geometric resolution. The construction of a decomposition mesh $\mathcal{K}_{D}$ is illustrated in Figure \ref{fig:HRmesh}.

\begin{figure}
    \centering
	 \begin{subfigure}[b]{0.3\linewidth}
	 \centering
    \includegraphics[width=\linewidth]{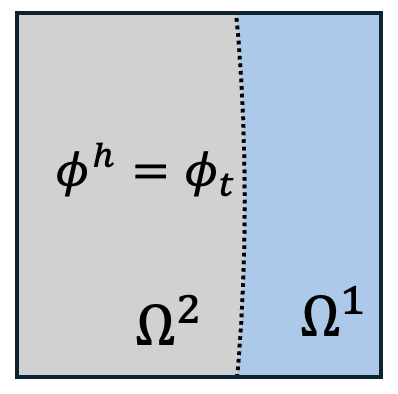}
    \caption{}
  \end{subfigure}
	 \begin{subfigure}[b]{0.3\linewidth}
	 \centering
    \includegraphics[width=\linewidth]{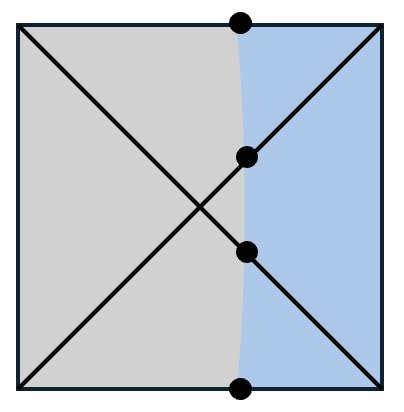}
    \caption{ }
  \end{subfigure}
  \begin{subfigure}[b]{0.3\linewidth}
	 \centering
    \includegraphics[width=\linewidth]{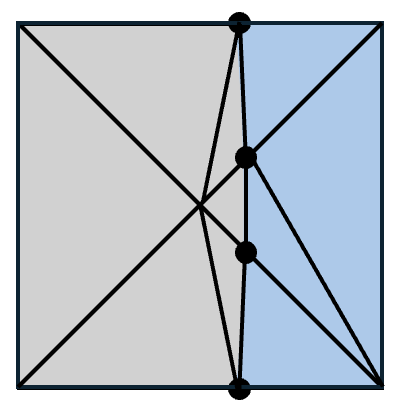}
    \caption{}
  \end{subfigure}
    \caption{ Foreground integration meshes are formed by triangulating cut  elements of the decomposition mesh $\mathcal{K}_{D}$. (a) A cell is intersected by the isocontour of the discretized level set function $\phi^h = \phi_t$ defining a material interface. (b) The cell is subdivided into triangular cells and isocontour-edge intersections (indicated by black cicles) computed. (c) Using the intersection points as nodal points, the cell is further subdivided. }
    \label{fig:triangulate}
\end{figure}

The decomposition mesh $\mathcal{K}_{D}$ is triangulated by first applying a pre-defined triangulation to the intersected background elements, as shown in Figure \ref{fig:triangulate} (a-b). This pre-defined triangulation forms 4 triangular elements in 2D domains or 24 hexahedral elements in 3D. Through root finding along elemental edges, the location of the isocontour is found, indicated by the black dots in Figure \ref{fig:triangulate} (b). A subdivision template is then applied to each intersected triangle/tetrahedron, as shown in Figure \ref{fig:triangulate} (c), to further subdivide the triangles/tetrahedrons into a set of triangles and tetrahedrons whose facets follow the LS isocontour.
The last step is repeated recursively for each individual LSF $\phi^h_i$ which enables sharp geometric corners and edges to be captured where multiple interfaces meet. This approach has been used previously by, e.g., \cite{soghrati_hierarchical_2014}.
The resulting approximation of the interface is piecewise linear and depends on the resolution of the decomposition mesh $\mathcal{K}_{D}$. 


Note that the sliver elements and poor aspect ratios in meshes produced by this method are still suitable for the interpolation of the background basis functions which is not bound by the typical mesh quality constraints of traditional FEM \cite{knupp_remarks_2007}. The resulting foreground mesh is of mixed element type (triangles and quadrilaterals in 2D, or hexahedrons and tetrahedrons in 3D) and contains hanging nodes. To accommodate the hanging nodes on the foreground mesh and to adequately interpolate the discontinuous enriched background basis functions described in Subsection \ref{subsec:Heaviside}, this method employs discontinuous Galerkin type elements for foreground function spaces. As the original THB-spline background basis maintains at minimum $C^{p-1}$ continuity within each material domain, the continuity of the interpolated basis is likewise $C^{p-1}$-continuous where interpolation is exact \cite{schillinger_lagrange_2016}. This work expands upon previous results using approximate extraction methods \cite{fromm_interpolation-based_2023}, where the constraints placed on the foreground Lagrange basis are lessened while numerical accuracy is maintained. 

The sliver elements resulting from this procedure do not themselves present problems with the interpolation-based workflow. However, due to the arbitrary location of material interfaces with respect to cell boundaries, it is possible for cells to be cut such that only a small portion of a basis function's support resides inside a given material domain. These small cell cuts result in sparsely supported basis functions, which can present issues with stability and linear conditioning. 

Numerous strategies exist to mitigate these issues and were recently reviewed in \cite{de_prenter_stability_2023}.
Strategies include basis function removal \cite{elfverson_cutiga_2018}, ghost stabilization \cite{burman_ghost_2010}, and basis function agglomeration \cite{badia_aggregated_2018} or extension \cite{burman_extension_2023}, and have the potential to be implemented within the presented interpolation-based framework. 
The benchmark problems presented in this work do not require special treatment of sparsely support basis functions and the authors leave the implementation of these stabilization strategies to future work. 


\section{Numerical Results}
\label{sec:NumResults}
The accuracy of the proposed method is demonstrated through the study of several benchmark problems. Problems were defined in the Python-based open-source FE code FEniCSx, using foreground meshes and extraction operators generated by the open-source XIGA code MORIS, available at \href{https://github.com/kkmaute/moris/}{github.com/kkmaute/moris} \cite{maute_moris_2023}. The source code with which the results were generated is also available at \href{https://github.com/jefromm/EXHUME_dolfinX}{github.com/jefromm/EXHUME\_dolfinX} \cite{fromm_jefrommexhume_dolfinx_2023}. 

\subsection{Resolving discontinuities in solution fields through Heaviside enrichment}
\label{subsecBar} 

\begin{figure*}
    \includegraphics[width=\textwidth]{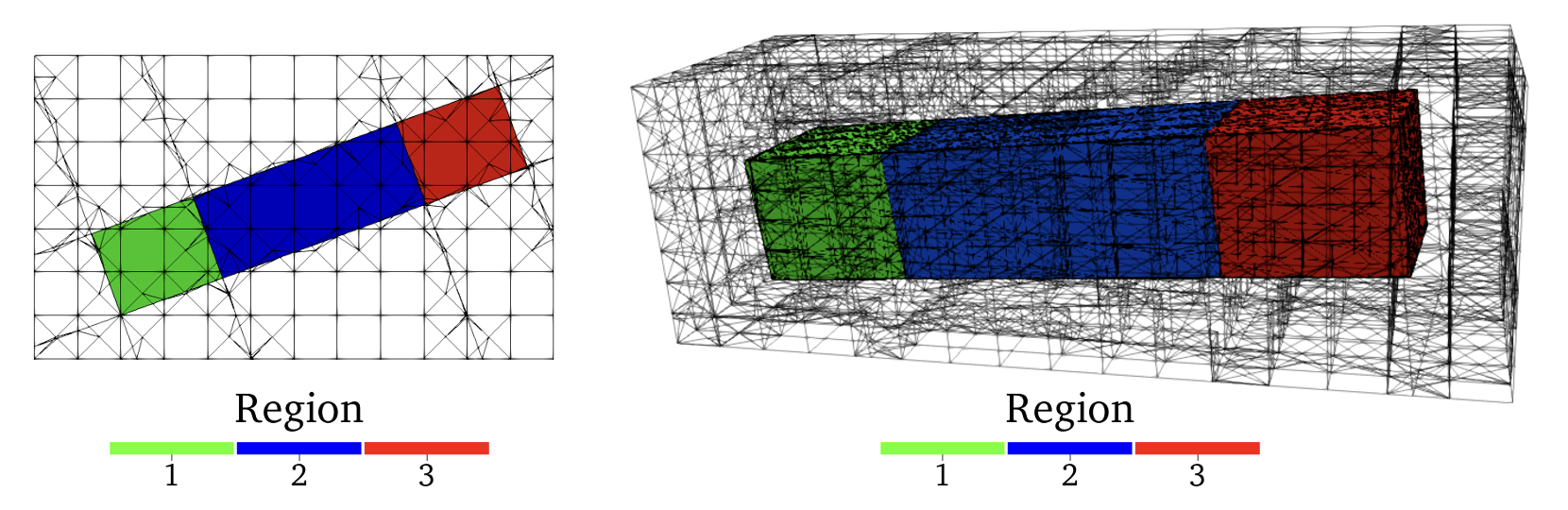}
    \caption{The geometric configurations for the 2D (left) and 3D (right) domains. A three-material beam is embedded in a structured background grid and rotated such that material interfaces do not align with element edges (in 2D) or facets (in 3D) . Elements intersected by the level set functions defining the beam geometry are triangulated to form a boundary-fitted foreground mesh.}
    \label{fig:barRegion}
\end{figure*}

\begin{figure*}
    \includegraphics[width=\textwidth]{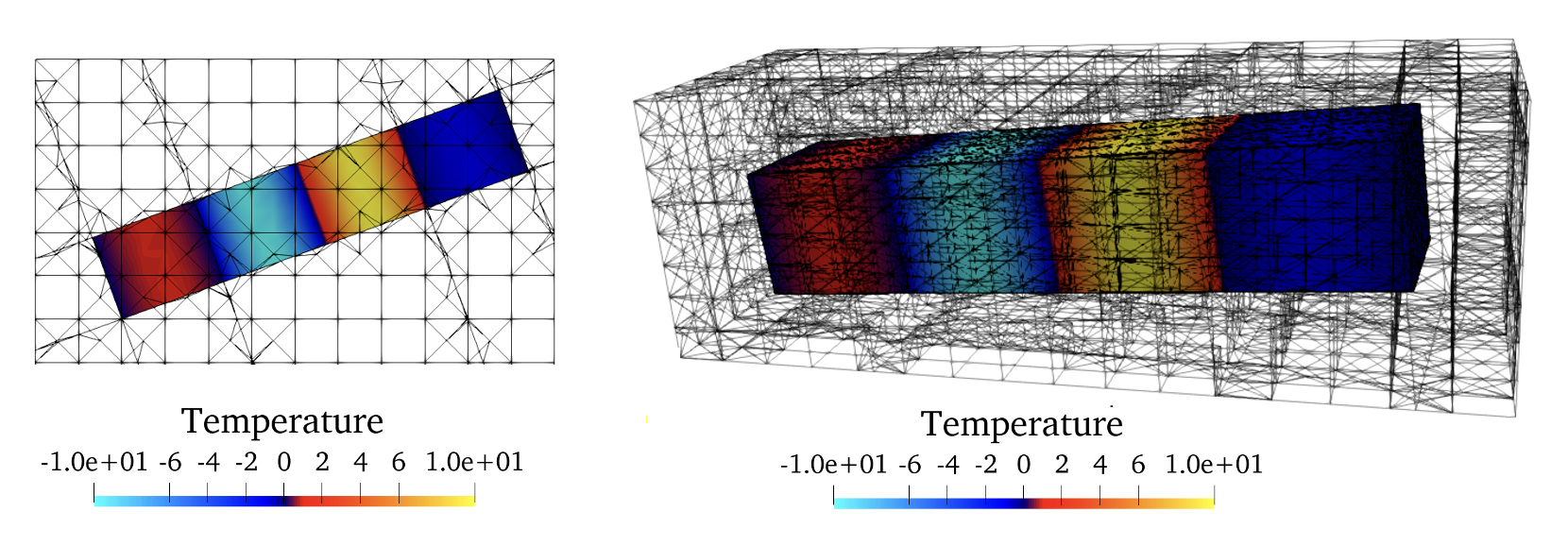}
    \caption{The temperature solution field plotted for both domain geometries, using a bi(tri)-quadratic B-spline basis interpolated with a bi(tri)-quadratic foreground Lagrange space. The solution is weakly discontinuous at material interfaces.}
    \label{fig:barTemp}
\end{figure*}

\begin{figure*}
    \centering
	 \begin{subfigure}[b]{\textwidth}
	 \centering
    \includegraphics[width=0.8\linewidth]{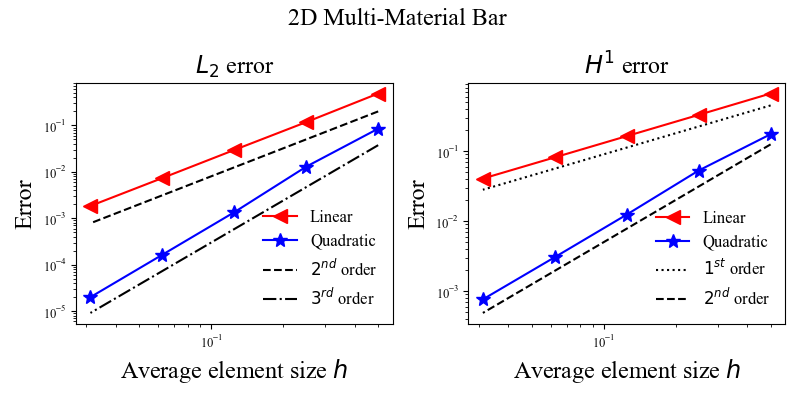}
  \end{subfigure}
	 \begin{subfigure}[b]{\textwidth}
	 \centering
    \includegraphics[width=0.8\linewidth]{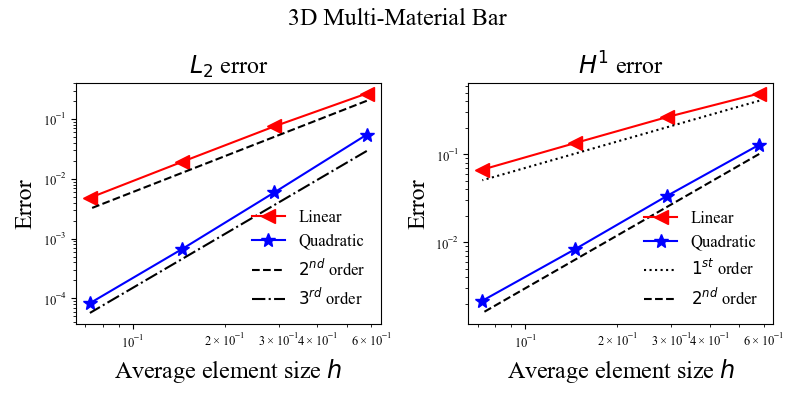}
  \end{subfigure}
    \caption{Ideal convergence rates are seen for both bi(tri)-linear and quadratic B-spline basis functions, which were interpolated with equal order foreground Lagrange function spaces. The convergence rates are ideal for both the 2D domain (top) and 3D domain (bottom).}
    \label{fig:barPlots}
\end{figure*}

In this Subsection, a multi-material beam undergoing a spatially varying heat load presents a weakly discontinuous temperature solution field. In this work, weak discontinuities refer to discontinuities in the gradients of solution fields. The solution is approximated with an interpolated Heaviside-enriched background basis, using the interpolation-based immersed boundary workflow. The error convergence results from this study demonstrate the accuracy of this method's enrichment scheme for multi-material problems. Beams in both 2D and 3D domains are considered. 


The beam is initially defined with corner coordinates $(0,0)$ and $(L,H)$ in 2D, and $(0,0,0)$ and $(L,H,H)$ in 3D, with $L=5$ and $H = 1$. The beam is divided into 3 sections, with interfaces at $x = L/4$ and $x= 3L/4$. Each section of the beam is assigned a thermal conductivity $\kappa_1 = 1.0$, $\kappa_2 = 0.1$, and $\kappa_3 = 1.0$.   To ensure the non-conformity of the material interfaces with respect to the background mesh facets, the 2D beam is rotated about the origin by angle $\phi = 20^\circ$, while the 3D beam is rotated about $y$- and $z$-axes by angles $\phi_y = -5^\circ$ and $\phi_z = 5^\circ$, respectively. The 2D beam is embedded into rectangle with corner coordinates (-1.0, -0.5) and (5,3) and the 3D beam is embedded into a rectangular prism with corner coordinates (-0.5, -0.25,-0.25) and (5.5,1.75,1.75). The geometric configurations are shown in Figure \ref{fig:barRegion}. Each edge (in 2D) or plane (in 3D) of the beam is defined by a level set function which extends beyond the beam domain in the mesh. The functions are extended through the entire mesh to fully resolve the corners (in 2D) or edges (in 3D) of the beam. 

Thermal diffusion is governed by the Poisson equation. The strong and weak forms of this problem are given by Eqs. \eqref{eq:poisson-strong} and \eqref{eq:heat-disc-generic} in Subsection \ref{subsec:ModelProblem}. 
The source term 
\begin{equation}
    q_{B} = -  \grad \cdot (\kappa \grad T_{ex}) 
\end{equation}
is constructed from the exact solution 
\begin{equation}
    T_{ex}(x') = \dfrac{1}{\kappa}\sin\left( \dfrac{4\pi}{L} x'\right).
\end{equation}
In 2D the beam-aligned coordinate is expressed in global coordinates as
\begin{align}
    x' = x \cos(\mbox{-}\phi) - y \sin(\mbox{-}\phi),
\end{align}
and in 3D
\begin{align}
    x' = z \sin(\mbox{-}\phi_y) + \big (x \cos(\mbox{-}\phi_z)\!-\!y \sin(\mbox{-}\phi_z)\big) \cos(\mbox{-}\phi_y),
\end{align}
where $\bm{x} = [x,y,z]$ are the mesh coordinates. The exact solution is imposed as Dirichlet boundary data on the ends of the bar, $x'=0$ and $x'=L$. 

For the 2D domains a suite of meshes with average background element sizes $h\in \cdot [1,0.5,0.25,0.125,0.0625]$ was used. In 3D, the background element sizes $h\in 1.16\cdot[1,0.5,0.25,0.125,0.0625,]$ were used. The foreground meshes were constructed with the workflow described in Subsection \ref{subsec:FGMesh}. Results from the 2D mesh with background element size $h = 0.5$ and from the 3D mesh with background element size $h = 0.25$ using bi(tri)-quadratic Lagrange foreground bases to interpolate a bi(tri)-quadratic enriched B-spline bases are shown in Figure \ref{fig:barTemp}.  Error convergence rates are plotted for both bi(tri)-linear and quadratic enriched background spline spaces in Figure \ref{fig:barPlots}, interpolated with equal-order foreground bases. Ideal error convergence rates validate the interpolation-based immersed boundary workflow for multi-material problems. 

\subsection{Approximating curved geometries through local foreground refinement }
\label{subsecEigenStrain}

\begin{figure*}
    \centering
\includegraphics[width=\linewidth]{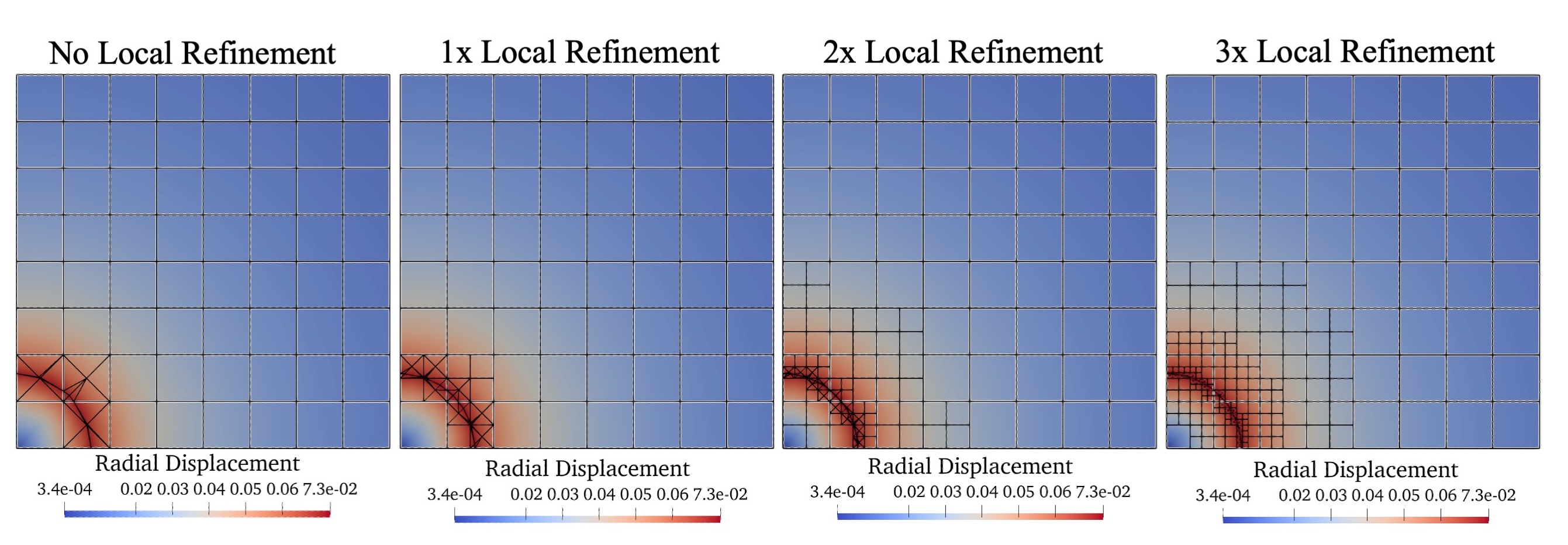}
    \caption{Approximated radial displacement of the eigenstrain problem. The background mesh is shown in white, while the foreground integration mesh is overload in black. The background mesh remains the same as local refinement is applied to the foreground for improved geometric resolution. The images correspond to the coarsest refinement level with background element size $h = 0.625$.}
    \label{fig:eigenFig}
\end{figure*}

\begin{figure*}
    \centering
\includegraphics[width=\linewidth]{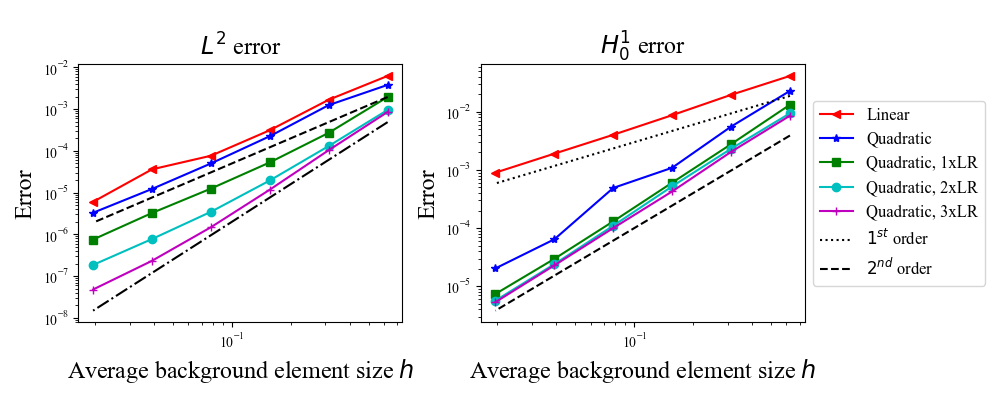}
    \caption{Error convergence data for the eigenstrain problem, illustrating the efficacy of foreground refinement. With no foreground refinement, the $L^2$ error convergence rate is limited by the geometric error. With sufficient foreground refinement (3x LR), the convergence rate approaches the ideal of 3. }
    \label{fig:eigenPlot}
\end{figure*}

A major challenge in the modeling of multi-material problems is the discretization of material interfaces. In this work, local refinement of the foreground mesh is performed to increase geometric resolution without affecting the number of degrees of freedom in the solution space. In this example the linear elastic behavior of an infinite plate with an embedded circular inclusion of radius $R=0.5$ is modeled, with local refinement employed to improve the approximation of the inclusion geometry.  

The inclusion is comprised of Material $1$ with Lam\'{e} constants $\lambda_1= 497.16$ and $\mu_1= 390.63$, while the exterior plate is made of Material $2$ with $\lambda_2 = 656.79$, and $\mu_2 = 338.35$.  A uniform isotropic eigenstrain of $\varepsilon_0=0.1$ is imposed on the inclusion. The weakly discontinuous analytic solution for the radial displacement is given in \cite{wang_homogenization_2003} as 
\begin{equation}
    u_r = \begin{cases}
    C_1 r, & r\leq R, \\
    C_1 \dfrac{R^2}{r}, & r\geq R, \\
    \end{cases}
\end{equation}
where 
\begin{equation}
    C_1 = \dfrac{(\lambda_1 + \mu_1)\varepsilon_0}{\lambda_1 + \mu_1 + \mu_2}.
\end{equation}

Exploiting symmetry, only the upper right quadrant of the plate is modeled as shown in Figure \ref{fig:eigenFig}. Symmetry conditions are enforced on the left and bottom edges of the domain, and the exact displacement is prescribed on the right and top edges. The solution domain is a $5\times 5$ square, with a quarter circle at the lower left corner.

The strong and discrete forms of the multi-material linear elasticity PDE are detailed in Appendix \ref{app:MMlinearElasticity}. For this example the mechanical strain is computed by 
\begin{equation}
    \bm{\varepsilon}_m (\bm{u}) = \begin{cases}
        \bm{\varepsilon}_u(\bm{u}) - \varepsilon_0 \bm{I} \text{, } &\bm{x} \in \Omega_1\\
        \bm{\varepsilon}_u(\bm{u}) \text{, } &\bm{x} \in \Omega_2\\\end{cases}
\end{equation} 
where the total strain $\bm{\varepsilon}_u(\bm{u}) = \frac{1}{2} ( \grad \bm{u} + (\grad \bm{u})^\text{T})$. 

The domain $\Omega$ is immersed into an axis-aligned $5 \times 5$ square on which the enriched background B-spline spaces are constructed. A suite of  background meshes with characteristic element lengths $h \in 0.625 \times [1, 0.5, 0.25, 0.125, 0.0625, 0.03125]$ are used to generate convergence data. 

The foreground meshes used here are locally refined about the material interfaces, as seen in Figure \ref{fig:eigenFig}. The background basis functions remain constant, meaning that there is no increase to the number of solution degrees of freedom with foreground refinement. 


The convergence plots in Figure \ref{fig:eigenPlot} show the expected error convergence rates for the bi-linear B-spline function space, while the $L^2$ error convergence of of the bi-quadratic function space is limited to the second order by the geometric error. This is confirmed by the reduction in error magnitude with local refined of the foreground mesh around the curved interface. With sufficient local foreground refinement the convergence rates approach the ideal rate. 

The foreground meshes generated using this octree local refinement strategy include hanging nodes, which are difficult for most commercial or open-source finite element software to handle. To accommodate these hanging nodes this interpolation-based workflow employs discontinuous foreground Lagrange polynomial basis functions. Classical discontinuous Galerkin methods require augmentation of the variational form to enforce higher levels of continuity at cell interfaces \cite{cockburn_discontinuous_2003}, but this augmentation is not done within this workflow. The results in this work indicate that an interpolated basis of lower continuity than required of a basis used for conforming finite element methods can be employed in this method, provided the background basis meets the continuity requirements. This attribute of interpolated-based methods was first exploited in \cite{fromm_interpolation-based_2023} to approximate fourth-order PDEs with a background basis of quadratic B-spline bases, which are $C^1$ continuous,  interpolated with a foreground basis of quadratic Lagrange polynomials, which are only $C^0$ continuous. 




\subsection{Image-based thermo-mechanical analysis of composite materials, utilizing multiple levels of local refinement}
\label{subsec:composite_image}

The capability of this method to tackle distinct discretization requirements of state variable fields within a multi-physics problem is illustrated here in a coupled thermo-elastic problem. This problem is posed on an alumina-epoxy composite sample undergoing simultaneous heating and loading conditions. Separate background discretizations are interpolated with a single foreground discretization for the two fields. 
\begin{figure}
\centering
	 \begin{subfigure}[b]{0.8\linewidth}
	 \centering
    \includegraphics[width=\linewidth]{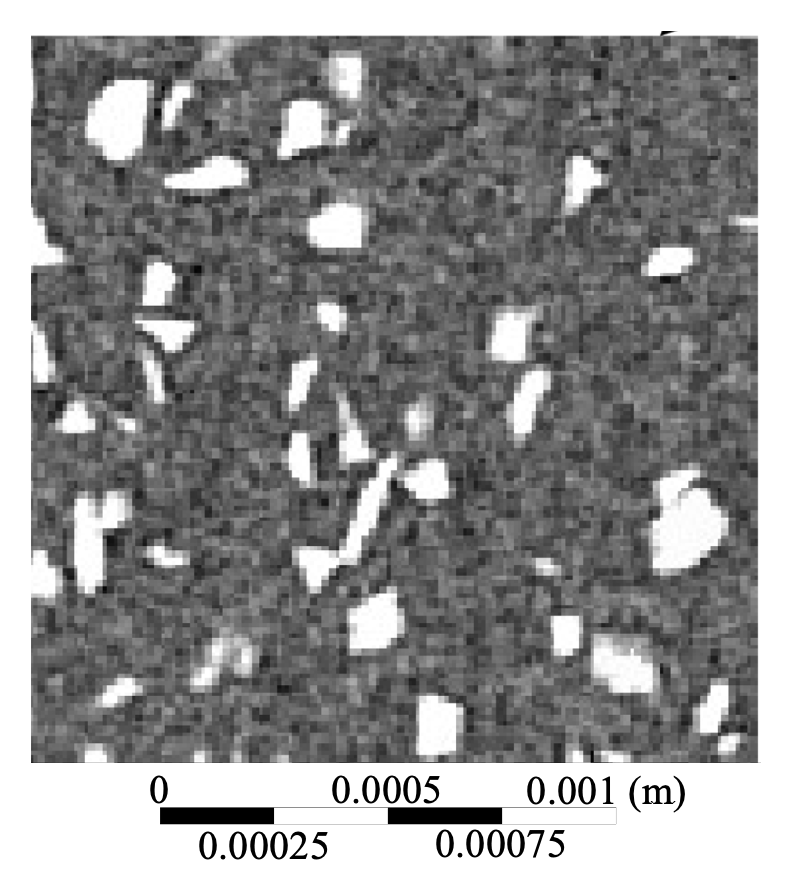} 
  \end{subfigure}
	 \begin{subfigure}[b]{0.8\linewidth}
	 \centering
    \includegraphics[width=\linewidth]{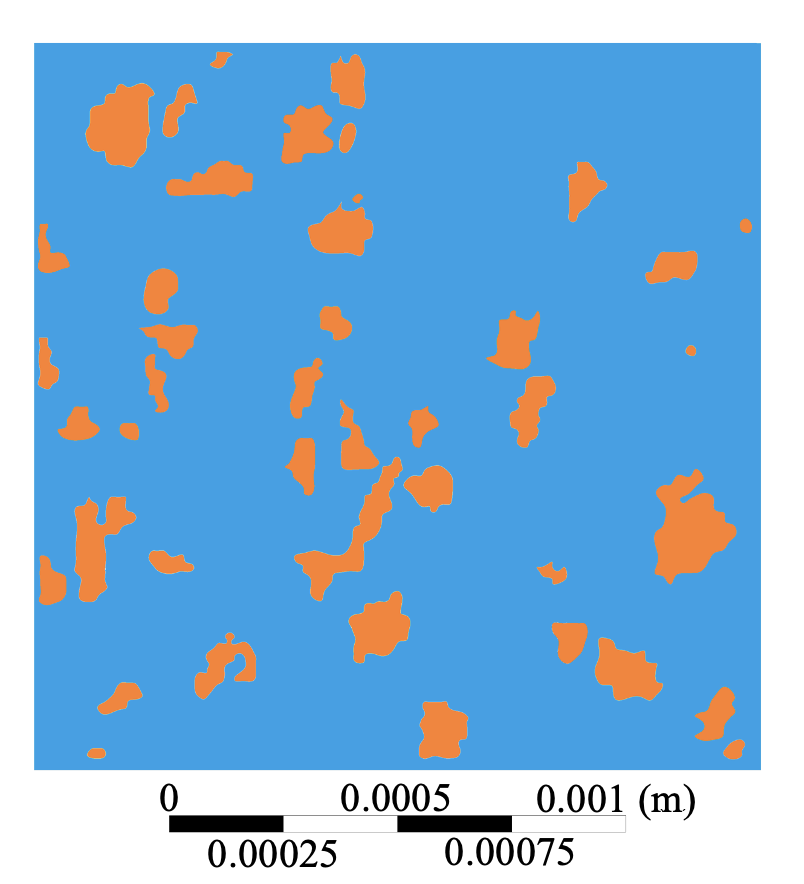} 
  \end{subfigure}
  \caption{Top:Micro-CT image of alumina-epoxy composite, where the white sections signify alumina particles and the grey is the surrounding epoxy.  Bottom: Smoothed image used to generate the LSF geometric description}
    \label{fig:microCT}
\end{figure}
\begin{figure*}
    \centering
	 \begin{subfigure}[b]{\linewidth}
	 \centering
    \includegraphics[width=\linewidth]{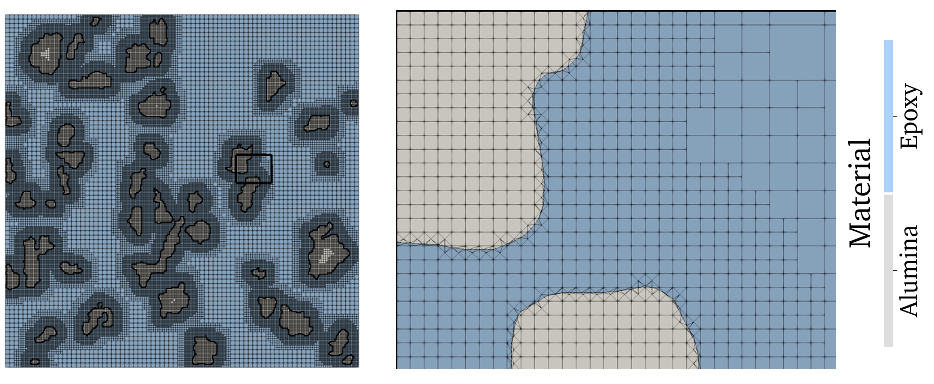} 
    \caption{The foreground integration mesh utilizes two levels of local refinement to resolve the geometry of the material interfaces, and contains 81,809 cells}
  \end{subfigure}
	 \begin{subfigure}[b]{1\linewidth}
	 \centering
    \includegraphics[width=\linewidth]{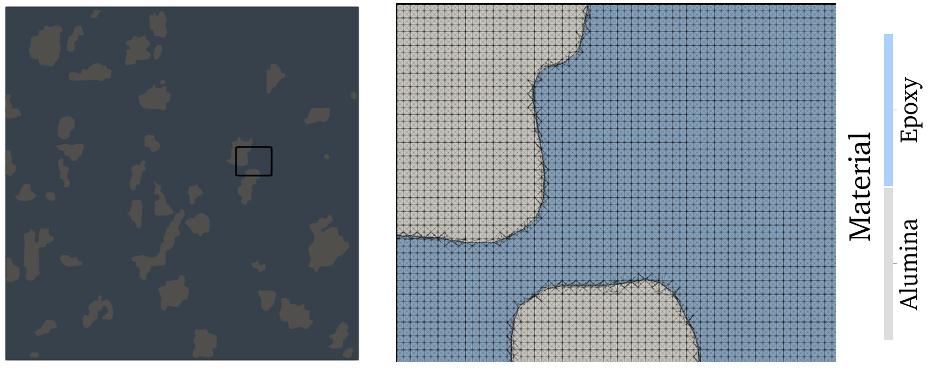} 
    \caption{The mesh used for the classical boundary-fitted FEM comparison was generated with 3 levels of uniform refinement. The mesh elements were uniformly triangulated and modified to improve mesh quality. The mesh contains 1,675,860 cells.}
  \end{subfigure}

    \caption{The whole domain is shown in the images on the left with the box indicating the region shown in the zoomed in view on the left.}
    \label{fig:compositeMeshes}
\end{figure*}

The geometry of the composite sample is taken from real micro-CT images converted to an implicit level-set description. The micro-CT image, taken from \cite{wang_support_2023} and shown in Figure \ref{fig:microCT}a, is made up of $200\times200$ pixels, with a pixel size of $8\mu$m, and the specimen is 1.6mm by 1.6mm. The epoxy is represented by the grey background while the alumina particles appear white. The image was then manually processed to generate the smoothed image shown in \ref{fig:microCT}b, which was then converted to the implicit level-set description. The following material properties, from \cite{master_bond_overview_2024,bauccio_ASM_1994}, were used: Poisson ratios $ \nu _{\text{Al}} = 0.23$ and $\nu _{\text{Ep}} = 0.358$, elastic moduli $ E _{\text{Al}} = 320e^9$Pa and $E _{\text{Ep}} = 3.66 e^9$Pa, thermal conductivities $\kappa _{\text{Al}} = 25.0$ W/mK and $ \kappa _{\text{Ep}} = 0.14$ W/mK, and thermal expansion coefficients $ \alpha _{\text{Al}} = 15e^{-6}$ 1/C$^o$ and $\alpha _{\text{Ep}} = 65e^{-6}$ 1/C$^o$. 

The process described in Subsection \ref{subsec:FGMesh} was used to construct a foreground integration mesh, beginning with a uniform 80 element by 80 element axis aligned decomposition mesh. Two levels of local refinement were applied about the material interfaces, and the refined decomposition mesh was triangulate. The resulting foreground mesh, shown in \ref{fig:compositeMeshes}a.  contains 81,809 cells. 

For comparison purposes,  a similar workflow was used to generate a boundary-fitted mesh for use in classical FEM. Classical FEM with FEniCSx requires a single element type mesh without hanging nodes. To avoid hanging nodes and to sufficiently resolve the gradients of the state variable fields the decomposition mesh was uniformly refined three times forming a 640 element by 640 element square mesh. The cut cells were triangulated to create a boundary-fitted mesh, and then the remaining quadrilateral elements were triangulated. The mesh was then modified with the software package Coreform Cubit 2023.11 to improve mesh quality metrics. The resulting mesh, shown in Figure \ref{fig:compositeMeshes}b, contains 1,675,860 cells. The method described here was used to ensure the level set descriptions of the material interfaces matched between this mesh and the foreground mesh used with the interpolation-based immersed boundary method.

\begin{figure*}
    \centering
    \includegraphics[width=\linewidth]{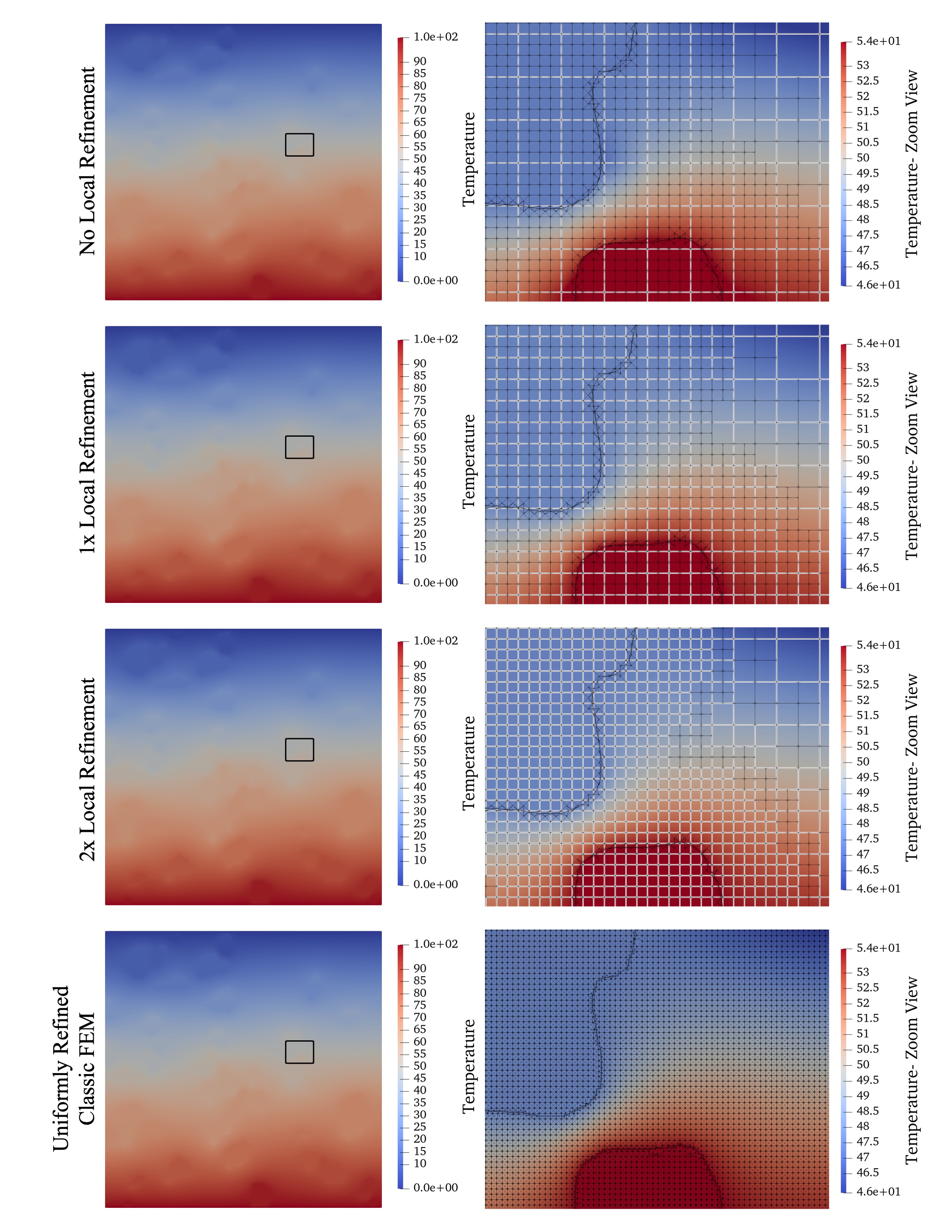}
    \caption{The temperature field results are compared for four different mesh configurations. The left shows the entire domain with the box indicating the region shown in the zoomed in view on the right.}
    \label{fig:compTemp}
\end{figure*}

\begin{figure*}
    \centering
    \includegraphics[width=\linewidth]{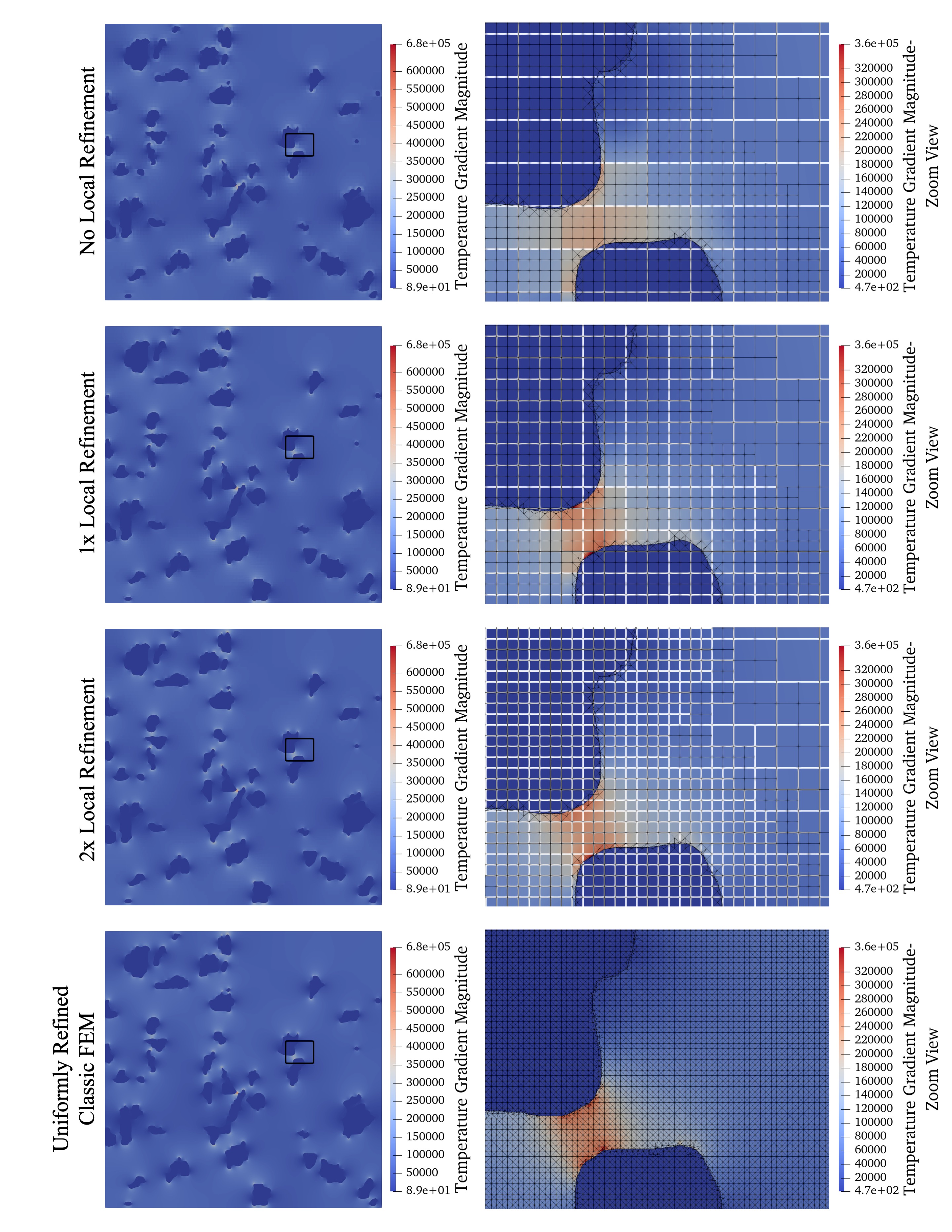}
    \caption{The the temperature gradient magnitude field is compared for four different mesh configurations. The left shows the entire domain with the box indicating the region shown in the zoomed in view on the right.}
    \label{fig:compTempGrad}
\end{figure*}

\begin{figure*}
    \centering
    \includegraphics[width=\linewidth]{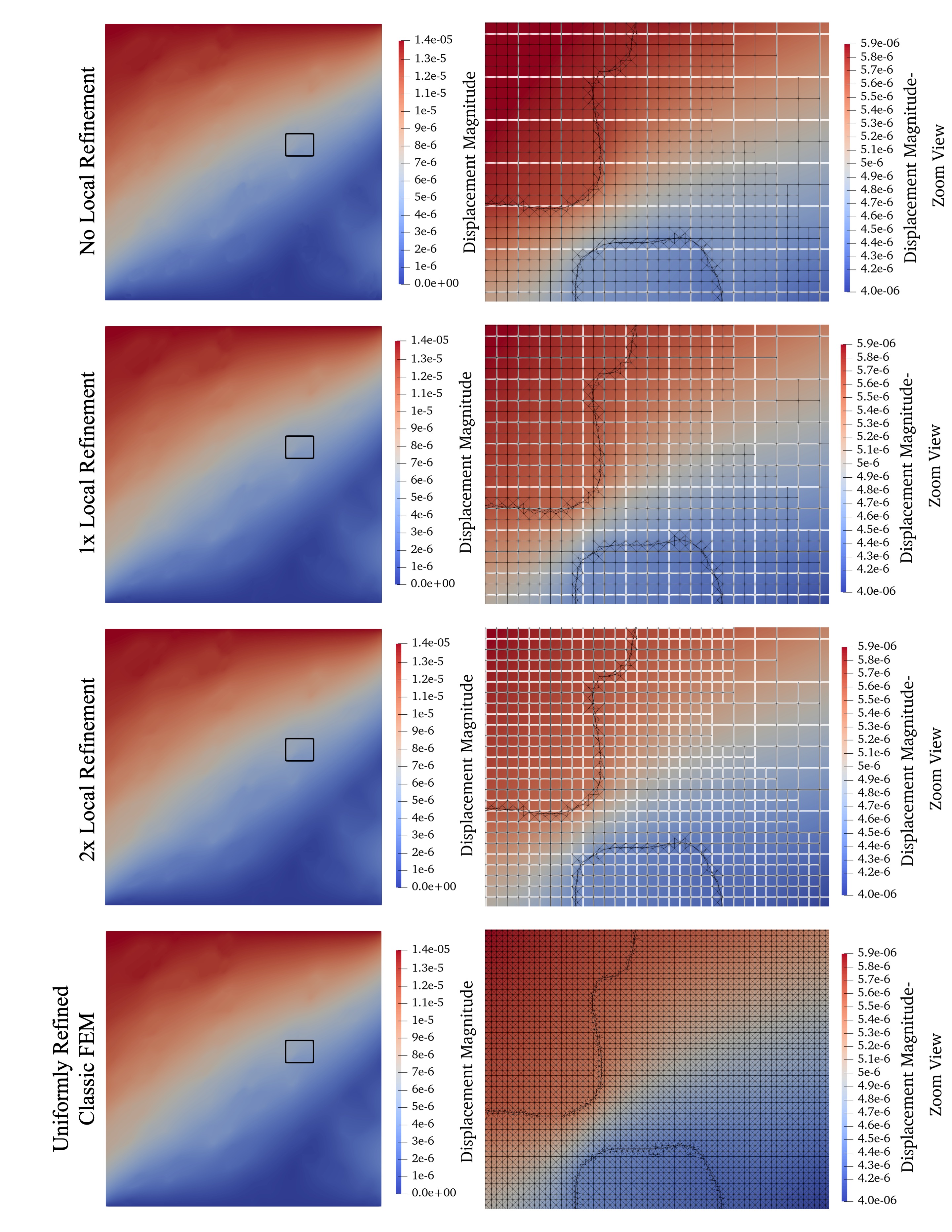}
    \caption{The displacement magnitude field results are compared for four different mesh configurations. The left shows the entire domain with the box indicating the region shown in the zoomed in view on the right.}
    \label{fig:compDisp}
\end{figure*}

\begin{figure*}
    \centering
    \includegraphics[width=\linewidth]{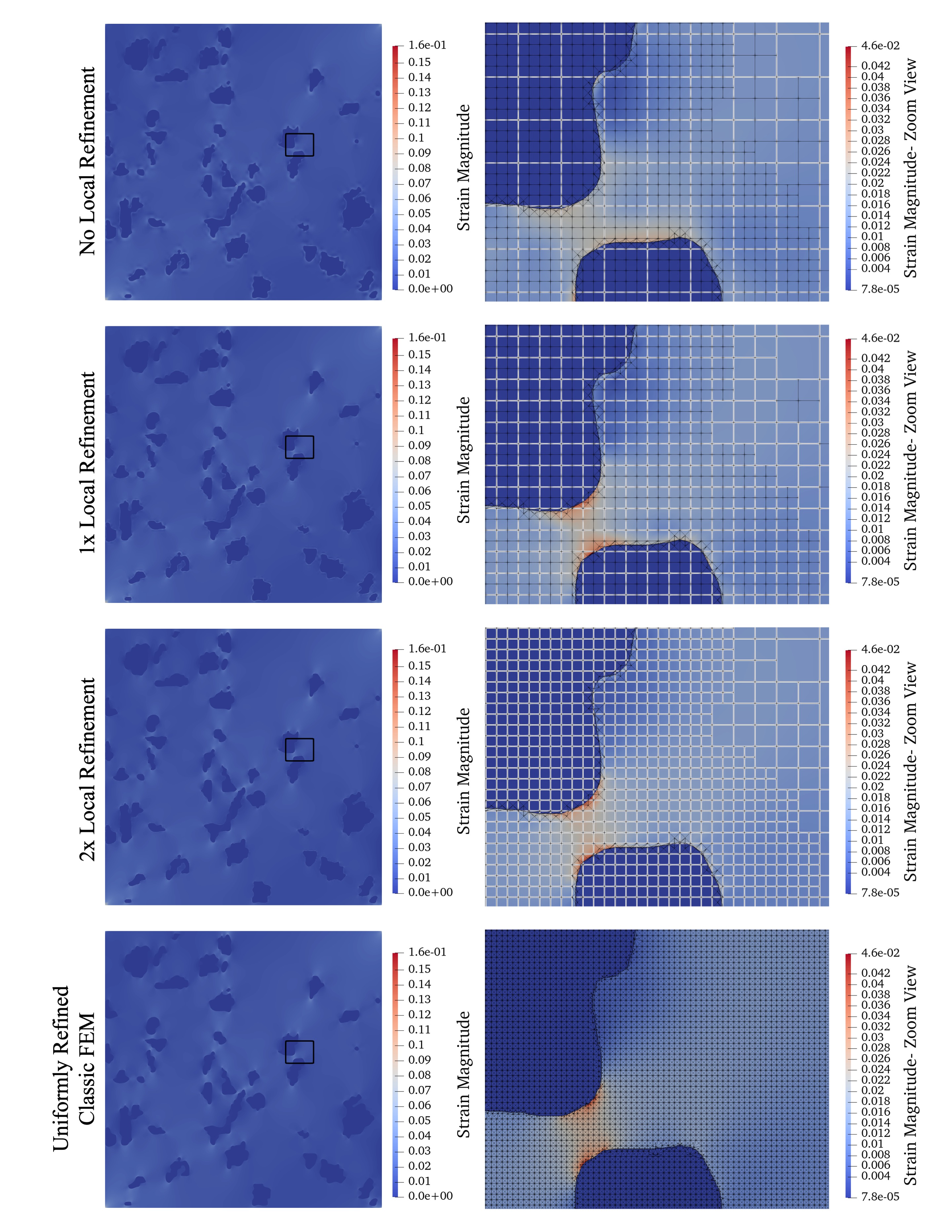}
    \caption{The mechanical strain magnitude field results are compared for four different mesh configurations. The left shows the entire domain with the box indicating the region shown in the zoomed in view on the right.}
    \label{fig:compStrain}
\end{figure*}

With this sample, a heated compression-shear test was simulated. The top and bottom displacements were imposed as $\bm{u}_{top} = [ -0.01, -0.01] $mm and $\bm{u}_{bottom} = [ 0, 0] $mm. The temperature at the top and bottom edges were specified as $T_{top} = 0^o$C and $T_{bottom}= 100^o$C. The environmental temperature was $T_0= 0^oC$. The sides were left both traction and heat flux free. 

The governing equations for heat conduction and linear elasticity were coupled via the constitutive law adding a thermal expansion component added to the total strain as in Eq. \eqref{eq:TEstrain}. Dirichlet boundary conditions for each field were imposed using Nitsche's method. Nitsche's terms enforcing continuity in $T$ and $\bm{u}$ were also imposed upon the allumina-epoxy interfaces. 

Bi-linear B-splines are used for the temperature field and bi-quadratic B-splines are used for each component of the displacement field, and both fields are interpolated using bi-quadratic Lagrange polynomials. For the boundary-fitted FEM comparison, Dirichlet boundary conditions were likewise enforced using Nitsches method, a bi-linear Lagrange basis was used for the temperature field, and a bi-quadratic Lagrange basis was used for the displacement. The results for the temperature field, the temperature gradient magnitude, the displacement magnitude, and the strain magnitude are shown in Figures \ref{fig:compTemp}, \ref{fig:compTempGrad}, \ref{fig:compDisp}, and \ref{fig:compStrain}, respectively. 
Within the figures, the foreground meshes are drawn in black while the background meshes are overlaid in white. The same foreground mesh, with two levels of local refinement to resolve the composite geometries, is used for each discretization.

\begin{table*}[h]
\centering
\caption{Degrees of freedom (DOFs) associated with the background function spaces at various levels of local refinement and with the uniformly refined FEM function spaces.}\label{tab:dofs}%
\begin{tabular}{@{}lrrrr@{}}
\toprule
       & No Local  & 1x Local  & 2x Local   & Uniformly Refined  \\
       &  Refinement &  Refinement &  Refinement & Boundary-Fitted Mesh  \\
\midrule
Bi-linear basis DOFs & 8,321 & 17,247 & 48,191  & 5,027,580   \\
Bi-quadratic basis DOFs & 9,381 & 22,166 & 60,729  & 10,055,160   \\
\hline
\end{tabular}
\end{table*}

Progressive levels of local refinement were applied to the background B-spline function spaces, as seen in the rows of Figures \ref{fig:compTemp}, \ref{fig:compTempGrad}, \ref{fig:compDisp}, and \ref{fig:compStrain}.  This local refinement was implemented with truncated hierarchically refined B-splines, as described in Subsection \ref{subsec:thb}. The regions of refinement with high-order splines are larger than those for the low-order to adequately support the truncated basis. This can be seen by comparing the white background meshes used for the temperature field depicted in Figures \ref{fig:compTemp} and  \ref{fig:compTempGrad} with those used for the displacement field in Figures \ref{fig:compDisp} and \ref{fig:compStrain}. The degrees of freedom associated with the various background discretizations are shown in Table \ref{tab:dofs}.

With two levels of local refinement, the results were in qualitative agreement with the ones of the uniformly refined classical FEM example with far fewer degrees of freedom. As previously noted, the uniformly refined meshes were initially generated with the same workflow used to create the boundary conforming foreground meshes, a `hands-off' method requiring only a bitmap image file. A more efficient mesh with fewer elements could be created for use in classical FEM but only with either significant user intervention or with sophisticated meshing software. Additionally, the workflow used here can easily be extended to 3D image stacks, which are not easily processed for classical FEM. 



\section{Conclusions}
\label{sec:Conclusions}

This work presents a new interpolation-based immersed boundary method for multi-material and multi-physics problems. This method employs  enriched truncated hierarchically refined B-spline background spaces and discontinuous hierarchically refined Lagrange integration spaces. 

Domain geometry and material interfaces are represented by level set functions, which can be generated by, for example, geometric primitives or from 2D or 3D images. The domain is embedded in a grid with an associated B-spline basis and the level set function is discretized and used to compute Heaviside indicator functions. Basis functions are inspected and individually enriched if their support is intersected by domain boundaries. This function-wise enrichment offers large savings in system degrees of freedom when compared to previous published on global enrichment strategies. 

Subdivision is used to generate a sequence of refined B-spline function spaces and associated tensor-product meshes from the original function space and mesh. In this work, the meshes are refined along domain interfaces. A hierarchically refined B-spline function space is defined recursively, along with its associated hierarchically refined mesh. This hierarchical function space is then truncated to form a partition of unity. 
The refined basis functions can then be enriched using the Heaviside enrichment functions when intersecting material interfaces. 

The enriched and refined background bases are interpolated by a discontinuous foreground basis. This foreground basis requires a boundary-fitted mesh, but is not subject to usual mesh conditioning constraints. The foreground meshes are thus constructed by the discretized level set function domain boundary descriptions and the hierarchically refined tensor-product meshes. Cells in the background mesh intersected by the domain and material interface geometry are triangulated to form a mixed-element type foreground mesh, which can be used by classical finite element codes to define a foreground basis. The level of refinement used to create the foreground mesh may exceed the level of refinement used for the background basis, allowing for greater geometric resolution without increasing the system's degrees of freedom associated with the background basis. 

The background basis is interpolated using extraction operators. Extraction operators are constructed by evaluating the background basis at the locations of the foreground basis nodes. Existing classical finite element codes assemble the linear systems using the Lagrange foreground basis. The linear system is then projected onto the interpolated basis with the extraction operators and the system is solved with the interpolated basis. Interpolation allows the enriched and refined B-splines basis to be utilized in traditional finite element codes without the addition of complex integration subroutines, broadening the applicability of this method.

Several benchmark problems validated this method. Interpolated enriched bases were used in both 2D and 3D for a multi-material thermal diffusion problem with exact geometric representation. The numerical approximations were compared with analytic solutions and errors were computed. The $L_2$ and $H_1$ errors converged at ideal rates with mesh refinement. Foreground only refinement was implemented for a multi-material linear elasticity PDE involving a circular domain. When compared to the analytic solution, errors from a non-refined foreground mesh converged ideally with a bi-linear basis, but the $L_2$ error convergence rate of the bi-quadratic basis was limited to $2^\text{nd}$ order due to geometric error in the discretization of the domain boundary. With sufficient foreground refinement, the ideal $3^\text{rd}$ order convergence for the bi-quadratic basis was observed. Lastly, micro-CT images were used to generate a geometric discretization of an alumina-epoxy composite sample. Thermo-elastity was simulated, coupling separately discretized temperature and displacement fields through a thermal expansion component. Unique truncated hierarchically refined B-spline bases were used for each field, bi-linear for temperature, and bi-quadratic for displacement. The results were in qualitative agreement with the ones of a uniformly refined traditional boundary conforming finite element simulation. 

In this work, the open-source code FEniCSx is used to demonstrate the efficacy of interpolation-based immersed boundary methods, and future work will expand implementation to other existing finite element codes. 
Within FEniCSx, additions to the interpolation-based immersed boundary workflow will be the implementation of stabilization techniques to address issues of linear conditioning.
\backmatter



\bmhead{Acknowledgements}
The first, second, third, and fourth authors acknowledge the support for this work from the National Science Foundation under Grants  2103939 and 2104106. The authors also acknowledge and thank Dr. David Kamensky for his assistance in developing the interpolation workflow. 






\begin{appendices}

\section{}\label{app:MMlinearElasticity}
\subsection*{Discretization of multi-material equations for linear elasticity}
Let a domain of $n$ material subdomains  $\overline{\Omega} = \overline{\Omega}^{1} \cup \overline{\Omega}^{2} \cup ... \cup \overline{\Omega}^{n}  \subset \mathbb{R}^d$ be the domain of interest, with varied material properties  $\mu$ and $\lambda$:
\begin{align}
\mu (\bm{x}) &= 
    \mu^m,  \bm{x} \in \Omega^m  
    \text{ and } \\
    \lambda (\bm{x}) =&
    \lambda^m,  \bm{x} \in \Omega^m .
\end{align}
Following the problem set up in Subsection \ref{subsecEigenStrain}, a source term $\bm{b}:\Omega\rightarrow\mathbb{R}^d$, a traction term $\overline{\bm{h}}:\partial \overline{\Omega}\rightarrow\mathbb{R}^d $ on $\Gamma_{\overline{\bm{h}}}\subset \partial \overline{\Omega}$ and Dirichlet boundary data $\overline{\bm{u}}:\partial \overline{\Omega}\rightarrow\mathbb{R}^d$ on $\Gamma_{\overline{\bm{u}}}\subset \partial \overline{\Omega}$ are ascribed. 

The strong form of this problem is then:  Find $\bm{u}: \Omega  \rightarrow \mathbb{R}^b$ such that $\forall$ $m \in \mathcal{M}$
\begin{equation}
\begin{split}
\label{eq:LEstrong}
    -\grad \bm{\cdot}\bm{\sigma}^m  = \bm{b}   ~~~ & \text{ in } \Omega^m \text{ ,} \\
    [\![ \bm{u}]\!] = 0    ~~~ &\text{ on all } \Gamma_{km} \\
    [\![ \bm{\sigma}]\!]\cdot \bm{n} = 0  ~~~ &  \text{ on all } \Gamma_{km} \\
    \bm{\sigma}^m\cdot\bm{n} = \overline{\bm{h}}\quad&\text{on}~\Gamma_{\overline{\bm{h}}}^m\text{ ,}\\
    \bm{u} = \overline{\bm{u}} \quad&\text{on}~\Gamma_{\overline{\bm{u}}}
\end{split}
\end{equation}
where $\Gamma_{km}= \overline{\Omega}^k \cap \overline{\Omega}^k \neq \emptyset$, $k \in \mathcal{M}$ and $ k \neq m$ are the material interfaces. $\Gamma_{\overline{\bm{h}}}^m = \Gamma_{\overline{\bm{h}}} \cup \partial \overline{\Omega}^m$, and $\Gamma_{\overline{\bm{u}}}^m = \Gamma_{\overline{\bm{u}}} \cup \partial \overline{\Omega}^m$ are intersections of the domain boundaries with the material subdomain boundaries, and $\bm{n}$ denotes the surface normal.  Here $[\![ \cdot]\!] = (\cdot)^{i} - (\cdot)^{j}$ is again the jump of a given quantity over the $\Gamma_{ij}$ interface, and $\bm{n}$ is the surface normal. For each material subregion the displacement is $\bm{u}^m =\bm{u}(\bm{x})$, $\bm{x} \in \Omega^m$.  The Cauchy stress tensor is defined as $\bm{\sigma}$
\begin{equation}
    \bm{\sigma}^m = \bm{C}^m : \bm{\varepsilon} =  2 \mu^m \bm{\varepsilon}^m + \lambda^m \text{tr}(\bm{\varepsilon}^m)\bm{I}
\end{equation}
in terms of the strain $\bm{\varepsilon}^m $, whose definition will vary depending on application. 

The computational domain $\Omega$ is embedded into a hierarchically refined background mesh $\mathcal{K}_u$, generated using a sequence of refined meshes $\mathcal{K}^l$ and subdomains $\Omega_{u}^l$, and associated with the  THB basis $\mathcal{T}_u = \{B^u_{i}\}$. Each component of the displacement is then discretized with the function space 
\begin{align}
    \mathcal{V}^h_u = \text{span}\{B^u_i| \text{ supp}(B^u_i) \cap \Omega \neq \emptyset\},
\end{align}
where $B^u_i \in \mathcal{T}_u$, the basis of enriched THB-splines. 

The discrete form is then: find $\bm{u}^h \in \bm{\mathcal{V}}_u^h= [\mathcal{V}_u^h, \mathcal{V}_u^h]$ such that all $\bm{v}^h \in \bm{\mathcal{V}}_u^h$
\begin{align}\label{eq:LE-disc-generic}
   \nonumber &\sum_{n=1}^m \left[\int_\Omega \bm{\sigma}(\bm{u}^h):\bm{\varepsilon}^u (\bm{v}^h) \,d\Omega \right]- \int_\Omega \bm{b} \cdot \bm{v}^h \,d\Omega \\
   &\quad - \int_{\Gamma_{\overline{\bm{h}}}}\overline{\bm{h}}\cdot \bm{v}^h\,d\Gamma\ = \mathcal{R}^{D}_u  + \mathcal{R}^{I}_u\text{ , }
\end{align}
where $\mathcal{R}^{D}_u$ and $\mathcal{R}^{I}_u$ are the Dirichlet and interface residuals, 
\begin{align}\label{eq:LEdirRes}
    \nonumber \mathcal{R}^{D}_u &=  \sum_{m=1}^n \Bigg[\mp \int_{\Gamma_{\overline{\bm{u}}}^m}(\bm{u}^h - \overline{\bm{u}})\cdot\bm{\sigma}(\bm{v}^h)\cdot\bm{n}\,d\Gamma\\
    \nonumber &\quad - \int_{\Gamma_{\overline{\bm{u}}}^m}\bm{v}^h\cdot\bm{\sigma}(\bm{u}^h)\cdot\bm{n}\,d\Gamma \\
    &\quad + \int_{\Gamma_{\overline{\bm{u}}}^m}\frac{\beta^{D}_u E }{h}(\bm{u}^h - \overline{\bm{u}}) \cdot\bm{v}^h\,d\Gamma\Bigg] \\
   \nonumber  \text{ and }&  \\
    \nonumber \mathcal{R}^{I}_u &=  \sum_{i=i}^n \sum_{j=1+i}^n \Bigg[ - \int_{\Gamma_{ij}}[ \![ \bm{u}^h]\!] \cdot(\{\bm{\sigma}(\bm{v}^h)\}\cdot\bm{n})\,d\Gamma \\
    \nonumber &\quad - \int_{\Gamma_{ij}}[ \![ \bm{v}^h]\!] \cdot(\{\bm{\sigma}(\bm{u}^h)\}\cdot\bm{n})\,d\Gamma \\
    &\quad + \int_{\Gamma_{ij}}\gamma^{ij}_u [ \![ \bm{u}^h]\!] \cdot[ \![ \bm{v}^h]\!]\,d\Gamma \Bigg]. 
\end{align}
As with the temperature Dirichlet residual, the first line opf Eq. \eqref{eq:LEdirRes} is either negative for symmetric Nitsche's method, which is used in this work, or positive for non-symmmetric Nitsche's method. Once again, $\{\cdot\} = w_{i}(\cdot)_{i} - w_{j}(\cdot)_{j}$ is the weighted average of a given quantity, with weights as defined in Eq. \eqref{eq:weights} using the elastic modulus as the material parameter $\omega$.  The penalty parameter $\gamma^{ij}_u $ is defined as 
\begin{equation}
    \gamma^{ij}_u = 2 \beta^{I}_u \dfrac{(h^{i})^{d_p-1} + (h^{j})^{d_p-1}}{(h^{i})^{d_p}  / E^{i} + (h^{j})^{d_p} / E^{j} }, 
\end{equation}
 where $\beta^{I}_u\geq 0$ is a user specified constant which controls the accuracy of the the interface condition, and $E(\bm{x}) = E^m$, $\bm{x} \in \Omega^m$. 

The linear elastic subproblem can be compactly written as the variational problem: Find $\bm{u}^h \in \bm{\mathcal{V}}^h_u = [\mathcal{V}^h_u, \mathcal{V}^h_u,]$ such that, $\forall \bm{v}^h\in \bm{\mathcal{V}}^h_u$
\begin{align}
    \label{eq:varU}
    a_u(\bm{u}^h, \bm{v}^h) = L_u(\bm{v}^h),
\end{align}
where $ a_u(\bm{u}, \bm{u})$ and $ L_u(\bm{v})$ can be derived from Eq. \eqref{eq:LE-disc-generic}. 

To approximate the solutions to multi-material linear elasticity problems, this workflow introduces an interpolated background function space 
\begin{align}
    \label{eq:interpolated_displ_space}
    \widehat{\mathcal{V}}^h_u &= \text{span}\{\widehat{B}^u_i \ \big| \  \text{supp}(\widehat{B}^u_i) \cap \Omega \neq \emptyset\},
\end{align}
where the interpolated background basis functions are defined 
\begin{align} 
    \widehat{B}^u_i &:= \sum_{j=1}^\nu M_{ij}^u N_j,
\end{align}
where 
\begin{align}
    M^u_{ij} := B^u_i(\bm{x}_j)
\end{align}
is the displacement component Lagrange extraction operator. $\{N_j\}_{j=1}^\nu$ is the basis of a Lagrange FE space with nodal points $\bm{x}_j$ such that $N_i(\bm{x}_j) = \delta_{ij}$. Here $\nu$ is the number of foreground basis functions. Note that for multi-physics problems, the same foreground space is used to interpolate the background bases for both the temperature and the displacements.

The approximation of each displacement component is then 
\begin{align}
    \label{eq:interpolationu}
    u_{k}^h  &= \sum_{i=1}^{n_u} \widehat{B}^u_i  d^{u_k}_i= \sum_{i=1}^{n_u} \sum_{j=1}^{\nu} M^u_{ij} N_j d^{u_k}_i, 
\end{align}
where $\{d^{u_k}_i\}_{i=1}^{n_{u}}$ are the unknown coefficients associated with each state variable field and $n_u$ is the number of basis functions in the interpolated background B-spline basis. $k \in \{ 1, ..., d_p\}$ are the indices associated with each of the displacement components, with $d_p$ denoting the physical dimension. The vector-valued approximation of displacement is defined as 
\begin{align}
    \label{eq:interpolationu}
    \bm{u}^h  &= \sum_{k=1}^{d_p} \sum_{i=1}^{n_u} \widehat{B}^u_i  d^{u_k}_i \bm{e}^k
\end{align}
where $\bm{e}^k$ are the directional unit vectors. For brevity in notation, new capital letter indices $I = d_p(i-1) +k$ and $J = d_p(j-1) +k$ are defined such that the vector value basis functions are 
\begin{align}
    \widehat{\bm{B}}^u_{I} = \widehat{B}^u_i \bm{e}^k \text{ and } \bm{N}_{J} = N^u_j \bm{e}^k 
\end{align}
and the approximation of displacement can be written 
\begin{align}
    \bm{u}^h  &= \sum_{I=1}^{(d_p \cdot n_u)} \bm{B}^u_{I}  d^{u}_{I} = \sum_{I=1}^{(d_p \cdot n_u)} \sum_{J=1}^{(d_p \cdot \nu)} M_{IJ}^{\bm{u}} \bm{N}_J   d^{u}_I
\end{align}
where $d^{u}_I = d^{u_k}_i$, and $M_{IJ}^{\bm{u}}$ are the components of the vector valued displacement field extraction operator
\begin{align}
    M^{\bm{u}}_{IJ} = \widehat{\bm{B}}^u _{I} \left(\bm{x}_J \right)
\end{align}
The variational form in Eq. \eqref{eq:varU} assembled using the interpolated basis forms the linear system 
\begin{align} \label{eq:linearSysU}
    \bm{K}^{\bm{v} \bm{v} } \bm{d}^{\bm{u}} = \bm{f}^{\bm{v} },
\end{align}
where 
\begin{align}
    K^{\bm{v} \bm{v} }_{IJ} &= a_{\bm{u}}(\widehat{\bm{B}}^u_I, \widehat{\bm{B}}^u_J), \text{ and } \\
    f^{\bm{v}}_{I} &= L_{\bm{u}}(\widehat{\bm{B}}^u_I,). 
\end{align}
Applying extraction, the linear system in Eq. \eqref{eq:linearSysU} is rewritten as
\begin{align} \label{eq:linearSysUEx}
    (\bm{M}^{\bm{v}} )^\text{T} \bm{A}^{\bm{v} \bm{v} } \bm{d}^{\bm{u}} =  (\bm{M}^{\bm{u}})^{\text{T}} \bm{b}^{u},
\end{align}
where $ \bm{A}^{\bm{v} \bm{v} } $ and $\bm{b}^{\bm{v} }$ are computed with the foreground basis 
\begin{align}
    A^{\bm{v} \bm{v} }_{IJ} &= a_{\bm{u}}(\bm{N}_I, \bm{N}_J) \text{ and} \\
    b^{\bm{v} }_{I} &= L_{\bm{u}}(\bm{N}_I).
\end{align}
The quantities $\bm{A}^{\bm{v} \bm{v} }$ and $\bm{b}^{\bm{v} }$ are evaluated and assembled with the boundary-fitted foreground mesh, following the workflow outlined in Subsection \ref{subsec:InterpolatedBF}. 

\end{appendices}


\bibliographystyle{bst/sn-basic}
\bibliography{MultiMatPaper}

\end{document}